\documentclass{m2an}
\usepackage{latexsym}
\usepackage{amssymb,amsmath,dsfont}
\usepackage{graphicx,float,color,subfigure,url}
\newtheorem{thm}{Theorem}[section]

\newtheorem{proposition}{Proposition}[section]

\newtheorem{rem}{Remark}[section]
\newcommand{\R}{\mathbb{R}}
\newcommand{\N}{\mathbb{N}}

\newcommand{\BB}{\boldsymbol{B}}
\newcommand{\vecdeux}[2]{\left(\begin{array}{c} #1\\#2 \end{array}\right)}
\newcommand{\WW}{\boldsymbol{W}}
\newcommand{\UU}{\boldsymbol{U}}
\newcommand{\VV}{\boldsymbol{V}}

\newcommand{\F}{\boldsymbol{F}}
\def\jump#1{\left[\left[ #1 \right]\right]}

\def\abs#1{\vert #1 \vert}
\def\div{{\rm div}}
\def\dsp{\displaystyle}

\newcommand{\resim}{F\textsc{igure}~}
\everymath{\displaystyle}
\begin{document}

\title{Air entrainment in transient flows in closed water pipes: a two-layer approach.}
\author{C. Bourdarias}
\address{Laboratoire de Math\'ematiques, UMR 5127 - CNRS and Universit\'e de Savoie, 73376 Le Bourget-du-Lac Cedex, France.\\
\email{Christian.Bourdarias@univ-savoie.fr, Stephane.Gerbi@univ-savoie.fr}
}
\author{M. Ersoy}
\address{BCAM--Basque Center for Applied Mathematics, Bizkaia Technology Park 500, 48160, Derio, Basque Country, Spain}
\secondaddress{\emph{Present address}: IMATH--Institut de Mathématiques de Toulon et du Var, Université du sud Toulon-Var,  Bâtiment U, 
BP 20132 - 83957 La Garde Cedex, France, 
\email{Mehmet.Ersoy@univ-tln.fr}
}
\author{S. Gerbi}
\sameaddress{1}


\begin{abstract}
In this paper, we first construct a model for  free surface flows that takes into account the air entrainment by a system of 
four partial differential equations. 
We derive it by taking averaged values of gas and fluid velocities on the cross surface flow in the Euler equations 
(incompressible for the fluid and compressible for the gas).  The obtained system is conditionally hyperbolic. 
Then, we propose a mathematical kinetic interpretation of this system to finally
construct a  two-layer kinetic scheme in which a special treatment for the ``missing'' boundary condition is performed. 
Several numerical tests on closed water pipes are performed and the impact of the loss of hyperbolicity is discussed and illustrated. 
Finally,  we make a numerical study of the order of the kinetic method in the case where the system is mainly
non hyperbolic. This provides a useful stability result when the spatial mesh size goes to zero.
\end{abstract}
\maketitle

{\bf Notations concerning geometrical quantities} \\
\begin{tabular}{p{0.1\linewidth}p{0.8\textwidth}}
$\theta(x)$ & angle of the inclination of the main pipe axis $z = Z(x)$ at position $x$ \tabularnewline
$Z(x)$ & slope \tabularnewline
$\Omega_{w}(t,x)$ &  cross-section area of the pipe orthogonal to the axis $z=Z(x)$ filled with water \tabularnewline
$\Omega_{a}(t,x)$ & cross-section area of the pipe orthogonal to the axis $z=Z(x)$ filled with air\tabularnewline
$\Omega(x)$ & $\Omega(x)= \Omega_w(t,x) \cup \Omega_a(t,x)$\tabularnewline
$S(x)$ &area of $\Omega(x)$\tabularnewline
$R(x)$& radius of the cross-section $\Omega(x)$\tabularnewline
$\sigma(x,z)$ & $\sigma(x,z)=\beta_r(x,z)-\beta_l(x,z)$ width of the cross-section $\Omega(x)$ at  altitude $z$ with $\beta_r(x,z)$ (resp. $\betaçl(x,z)$) is the right (resp. left) 
boundary point at altitude $z$.
\end{tabular}

{\bf Notations concerning the water model} \\
\begin{tabular}{p{0.1\linewidth}p{0.8\textwidth}}
$\rho_0$& density of the water at atmospheric pressure $p_0$ \tabularnewline
$H_w(x,z)$ & water height at section $\Omega_{w}(t,x)$ \tabularnewline
$\dsp A(t,x)$ &  the wet area\tabularnewline
$u(t,x)$ & velocity\tabularnewline
$Q(t,x)$& $Q(t,x) = A(t,x) u(t,x)$ is the discharge\tabularnewline
$c_w$ & air sound speed\tabularnewline
\end{tabular}

{\bf Notations concerning the air model} \\
\begin{tabular}{p{0.1\linewidth}p{0.8\textwidth}}
$\gamma$ & adiabatic index $\gamma$, set to $\dsp7/5$ for a diatomic gas\tabularnewline
$\mathcal{A}(t,x)$ & $\mathcal{A}(t,x) = S(x)-A(t,x)$ air cross-section area\tabularnewline
$\rho_a$ & constant reference density\tabularnewline
$p_a$ & constant pressure reference $\Omega(x)$ at  altitude     $z$\tabularnewline
$\rho(t,x,y,z)$ & density of the air at the current pressure\tabularnewline
 $\dsp\overline{\rho}(t,x)$ & $\dsp\overline{\rho}(t,x)= \frac{1}{\mathcal{A}(t,x)}\int_{\Omega_a(x)} \rho(t,x,y,z)\,dy\,dz$ 
 is the mean value of $\rho$ over $\Omega_a(t,x)$\tabularnewline
$\rho_0$ & density of the water at atmospheric pressure $p_0$\tabularnewline
$H_a(x,z)$ & air height at section $\Omega_{a}(t,x)$ \tabularnewline
$M(t,x)$& $M(t,x) = \frac{\overline{\rho}(t,x)}{\rho_0} \mathcal{A}(t,x)$  pseudo area associated to the air layer\tabularnewline
$v(t,x)$ & mean air speed \tabularnewline
$D(t,x)$ & $D(t,x)=M(t,x)v(t,x)$ pseudo air discharge\tabularnewline
$c_a$ & air sound speed\tabularnewline
\end{tabular}

Bold characters are used for vectors and $\alpha$ for both layers when no ambiguities are possible, otherwise $\alpha=a$ for the air layer or $\alpha=w$ for the water layer.

\section{Introduction}
The hydraulic transients of flow in pipes have been investigated by Hamam and McCorquodale \cite{HM82}, Song et al. \cite{S76,S77,SCL83}, 
Wylie and Streeter \cite{WS93} and others.  
However, the air flow in a pipeline and  the associated air pressure surge was not studied extensively. There is a large literature in which
several mathematical models and numerical schemes are presented. Difficulties such as the ill-posedness, the presence 
of discontinuous fluxes and also the loss of hyperbolicity are studied. This framework enters also in the class of the two-layer shallow water equations. 
Let us make of a short review of the models that take  into account the air and the problem encountered. 

The two component gas-liquid mixture flows occur in piping systems in several industrial areas such  as
nuclear power plants, petroleum industries, geothermal power plants, pumping stations and sewage pipelines.
Gas may be entrained in other liquid-carrying pipelines due to cavitation, \cite{WS79} or gas release from solution due to a drop
or increase in pressure. 

Unlike in a pure liquid in which the pressure wave velocity is constant, the wave velocity in a gas-liquid mixture varies with the pressure. Thus the main
coefficients in the conservation equation of the momentum are pressure dependent and consequently the analysis of transients in the two-component flows 
is more complex.
 
The most widely used analytical models for two-component fluid transients are the homogeneous model \cite{CBMN90}, the drift flux model 
 \cite{FH99,HI03,I75}, 
and the separated flow model \cite{TP97,CPT01}.
 
In the homogeneous model, the two phases are treated as a single pseudo-fluid with averaged properties \cite{CBMN90,WS93}:  
there is no relative motion or slip between the two phases.
The governing equations are the same as the one for a single phase flow. However the inertial and gravitational effects can play a more important role
than the relative velocity between the air phase and the liquid one: the interaction between the two phases has to be taken into account.

In the drift-flux model, \cite{FH99,HI03,I75}, the velocity fields are expressed in terms of the mixture center-of-mass velocity and the drift velocity of the air phase, which is 
the air velocity with respect to the volume center of the mixture. 
The effects of thermal non-equilibrium are accommodated in the drift-flux model by a constitutive equation for phase change that specifies the rate of mass transfer per unit volume. 
Since the rates of mass and momentum transfer at the interfaces depend on the structure of two-phase flows, these constitutive equations for the drift velocity 
and the vapor generation are functions of flow regimes.

Two-phase flows always involve some relative motion of one phase with respect to the other; therefore, a two-phase flow problem should be formulated in terms of two velocity fields. A general transient two-phase flow problem can be formulated by using a two-fluid model. The separated flow model considers the two phases individually, interacting with each other. Generally, this model will be written as a system 
of six partial differential equations  that represent  the balance of mass, momentum and energy \cite{CPT01,TP97}. 

A two-phase fluid flow composed of liquid droplets, assumed incompressible, suspended in  a compressible gas 
are governed by an Euler system where each phase are seen as a fluid described by macroscopic quantities is considered in 
\cite{Sainsaulieu93}. In this framework, the mass density of the liquid is very large in comparison with the gas mass and have application 
for cryogenic purposes.  The Euler system of six equations, written in a nonconservative form, 
is hyperbolic for some region depending on the volume fraction,  
the speed of each phases and the gas sound speed. Sainsaulieu \cite{Sainsaulieu95}  introduces also 
a finite volume approximation based 
on an approximate Roe-type Riemann solver which provides good results  with the exact solution in one dimension slab geometry for a two phase 
fluid flows  model conditionally hyperbolic.

In \cite{SW84}, the link between two-phase flows and flows in pipes has been studied and the resulting hydrodynamic equations
may also possess complex eigenvalues. As a consequence the initial-value problems based on such equations are ill-posed and
any consistent finite difference numerical scheme for these equations is unconditionally unstable, i.e., for
any constant ratio $\Delta t/\Delta x$, geometrically growing instabilities will always appear if $\Delta x$ 
is sufficiently small. 

Difficulties encountered in this paper enter also  in the framework of the so-called two-layer shallow water equations. 
While nonconservative term implies the lack of appropriate Rankine-Hugoniot relation, the conditionally hyperbolic condition may 
introduce complex eigenvalues as in the two-phase flows.  Moreover, the system is only hyperbolic for small relative 
speed of each layers for the mode called external wave motion while the internal is specific to these complex eigenvalues.  
As pointed out by several authors, the violation of such condition  
is linked to  well-known Kelvin-Helmholtz instability, i.e. corresponding to the strong velocity shear between layers, for which 
the vertically averaged two-layer system is not \emph{a priori} suitable (see  for instance \cite{AK09}). 

Several attempts have been proposed to solve these problems. We may cite the numerical scheme proposed by Abgrall and Karni \cite{AK09} which 
deals with nonconservative and non hyperbolic system by using relaxation approach which offers a greater decoupling and provides an access to the eigenstructure.
We refer to \cite{ParesRev06}  for discussion of numerical schemes and to \cite{B04} 
for the recent progress on hyperbolic systems with source terms.

Finally, although questions about the ill-posedness equations remain unsettled, such \emph{a priori} ill-posed equations are
widely used in practice such as for multiphase-flow computations,  backward heat conduction and porous media flows.

In the present paper, we consider the air entrainment appearing in the transient flow in closed pipes filled by  the free surface flow and air flow.
We derived the present model in a very closed manner as the derivation of 
Saint-Venant equations for free surface open channel flows: we consider that the liquid is inviscid and incompressible (thus Euler equations are used)
and its pressure is only due to the gravitational effect. 
It is the well known hydrostatic pressure law. A kinematic law at the free surface and a no-leak condition are used to complete the Euler system.
The air phase is supposed to be compressible and isothermal: thus we can use an equation of state
of the form $p = k \rho^\gamma$ where $\gamma$ represents the adiabatic exponent. 
To connect the two phases, we write the continuity of the normal stress tensor at the free surface separating
the two phases: the hydrostatic  pressure law for the fluid is thus coupled to the pressure law governing the air flow at the free surface. 
Then by taking averaged values over cross-sections of the main flow axis, 
we obtain a system of four partial differential equations which is conditionally hyperbolic: 
we  call it the two-layer water/air model. In Theorem \ref{BilayerModelThm}, we state  
its main property: the existence of a total energy. 
The derivation of the proposed  model and the analysis of the mathematical properties are presented in Section \ref{Model} 
 as well as a careful study of the eigenvalues of the air-layer model:  the problem is hyperbolic for small relative speed but also for large relative speed which was already pointed out
in the pioneer work of \cite{Ovsjannikov79} and recalled later on by Barros and Choi \cite{BC08} and the reference therein.

As in \cite{PS01,BEG11_1,BEG11_2}, we choose to interpret this model as a kinetic model since 
the computation of eigenvalues are not required: we present the mathematical kinetic formulation 
of the system of four partial differential equations in Section \ref{Kinetic}. 

This mathematical formulation allows us to introduce the two-layer kinetic scheme presented in  Section \ref{KineticScheme} which allows 
to reproduce fairly physical solutions even if the system loses its hyperbolicity.  
To deal with real boundary conditions, we also introduce in Section \ref{SectionBoundaryTreatment} the ``decoupled boundary conditions'' based on 
the solution of the Riemann problem for the decoupled system for which the eigenstructure is completely well-defined.

Finally, we propose in Section \ref{Numeric} some numerical simulations  to focus on the role of the air entrainment on the free surface fluid flow 
and the behavior of the numerical solutions provided for internal wave motion, i.e. when the system is not hyperbolic. 
More precisely, even if the presented two-layer system is not hyperbolic,  the numerical scheme seems to be stable and there are no geometric instabilities 
for a large series of numerical tests. 
This surprising stability is certainly due to the fact that the energy is conserved as stated in Theorem \ref{BilayerModelThm}. 
The scheme is then stable with respect to smaller $\Delta x$ and produces, even is the decay rate 
is worse than the case where the air is not taken into account, usable numerical results.

\section{A two-layer water/air model}\label{Model}
We will consider  throughout this paper a transient flow in a closed pipe that is not completely filled of water nor of air. 
The flow is composed by two air and water  layers separated by a moving interface. We assume that  the two layers are immiscible.
We will consider that the fluid is perfect and incompressible whereas the air is a perfect compressible gas: 
we neglect here the thermodynamic relations such as the air-water interactions  as 
condensation/evaporation since we suppose that the air layer follows an isothermal process and is assumed to be isentropic.

\begin{figure}[H]
\subfigure[Geometric characteristics of the domain \label{OxOz}]{\includegraphics[height = 6cm]{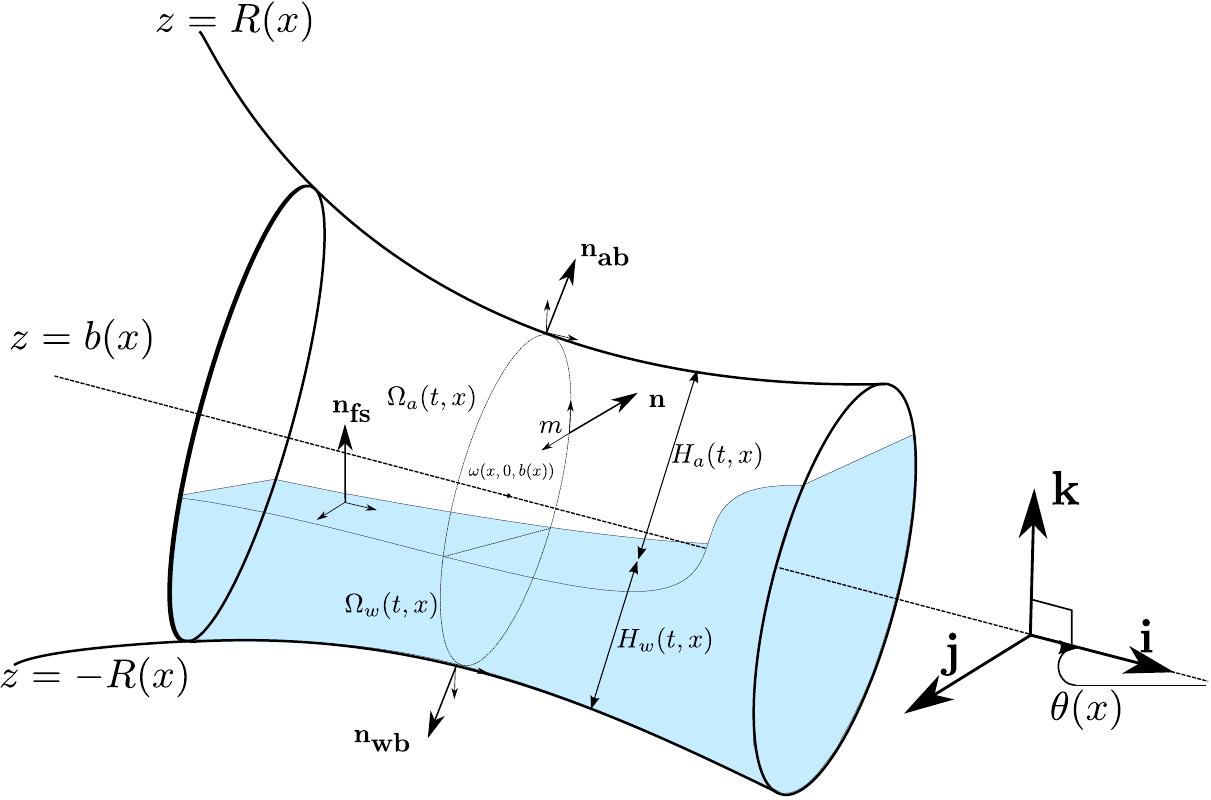}}
\subfigure[Cross-section of the domain \label{OyOz}]{\includegraphics[height = 5cm]{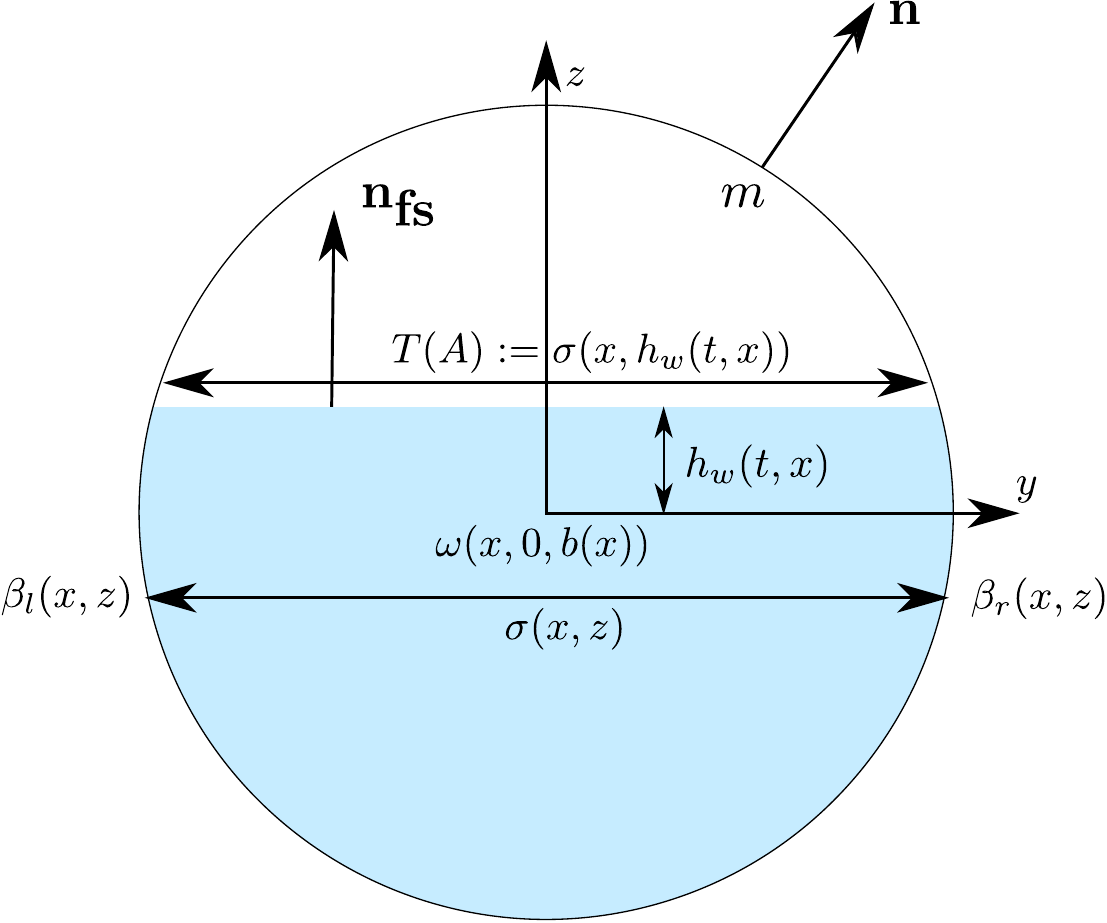}}
\caption{Geometry of the pipe: notations and settings\label{Geom}}
\end{figure}

The derivation of the two-layer model is then performed with the use of the 3D Euler (incompressible for the fluid, compressible for the gas) equations
where the no-leak condition is assumed for the fluid as well as for the gas at the pipe boundary.
The continuity of the pressure at the interface separating each layers is assumed.  
\begin{rem}
Throughout the paper, we use the following notations: 
we will denote by the index $w$ the water layer and by the index $a$ the air layer. 
We also use the index $\alpha$ to write indifferently the liquid or gas layer.
\end{rem}
The pipe is assumed to have a symmetry axis $\mathcal{C}$, says circular one, 
which is represented by a straight line with a constant angle $\theta$. 
$(O,\mathbf{i},\mathbf{j},\mathbf{k})$ is a reference frame attached to 
this axis with $\mathbf{k}$ orthogonal to $\mathbf{i}$ in the vertical plane containing 
$\mathcal{C}$,  see \resim \ref{OxOz}.
We denote by $\Omega_{t,a}$ and $\Omega_{t,w}$ the gas and fluid domain occupied at time $t$ as represented 
on \resim\ref{Geom}.
Then, at each point $\omega (x)$, we define the water section  $\Omega_{w}(t,x)$ by: 
$$\Omega_{w}(t,x)=\left\{(y,z)\in\R^2; z \in [-R(x),-R(x)+H_{w}(t,x)],\, y \in [\beta_l(x,z),\beta_r(x,z)]\right\}$$ 
and the air one by:
$$\Omega_{a}(t,x)=\left\{(y,z)\in\R^2; z \in [-R(x)+H_{w}(t,x),R(x)],\, y \in [\beta_l(x,z),\beta_r(x,z)]\right\}$$
where $R(x)$ stands for the radius, 
$H_{w}(t,x)$ the water height at section $\Omega_{w}(t,x)$, $H_{a}(t,x)$ the air height at section $\Omega_{a}(t,x)$.
$\beta_l(x,z)$, $\beta_r(x,z)$ are respectively the left and right $y$ values of the  boundary points of the domain at altitude $-R(x)<z<R(x)$ (see \resim \ref{OxOz} and  \resim \ref{OyOz}). 

We have then the first natural  coupling:
\begin{equation}\label{AirWaterRelation}
H_{w}(t,x)+H_{a}(t,x)=2R(x)\,.
\end{equation}
In what follows, we will denote $h_{w}(t,x)=-R(x)+H_{w}(t,x)$ the algebraic water height.
\begin{rem}\label{remtheta1}
Let us emphasize that the assumption  $\theta = cst$ is made for the sake of simplicity since the variable case does not introduce 
novelty with respect to the reference therein. 
Nevertheless, let us note that the curved ground case (see, for instance, 
\cite{BEG2012,SavageHutter} and \cite[Chapter 1]{TheseErsoy}) 
is adapted by introducing the  assumption on the curvature radius $\mathcal{R}(x)$, i.e. $\frac{\mathcal{R}(x)}{R(x)}>1$ where $R(x)$ 
is the radius of the pipe at $x$ along the curvilinear frame. 
This hypothesis ensure that the transformation of the initial Euler equations written in a given frame to the curvilinear one is a
 diffeomorphism  
and avoid the situation as illustrated in \resim \ref{FigHypothesisH}. We refer to \cite{BEG2012} for a detailed 
derivation of it.
\begin{figure}[H]
\includegraphics[height=5cm]{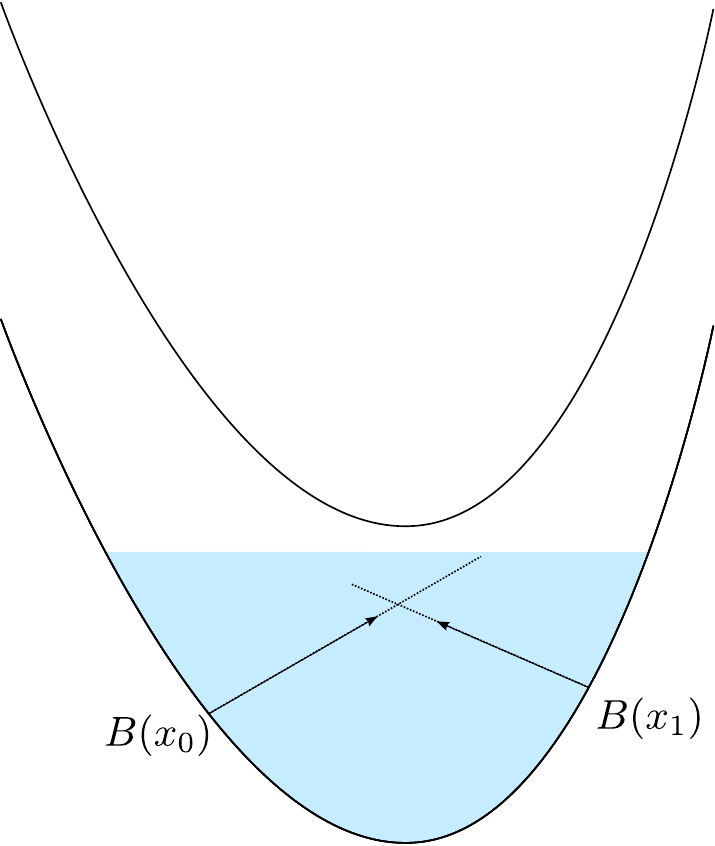}
\caption{Geometric restriction on the pipe for curved ground\label{FigHypothesisH}}
\end{figure}
\end{rem}

\subsection{The fluid free surface model}\label{FluidModel}
The shallow water equations for free surface flows are obtained from  the incompressible Euler equations \eqref{EulerModelFluid} which writes in cartesian coordinates $(x,y,z)$:
\begin{equation}\label{EulerModelFluid}
\begin{array}{llll}
\div (\mathbf{U_{w}}) &=& 0, & \textrm{ on } \R \times\Omega_{t,w} \\
\partial_t (\mathbf{U_{w}}) + \div (\mathbf{U_{w}} \otimes \mathbf{U_{w}}) + \nabla P_{w} &=&\mathbf{F}, & \textrm{ on }\R\times\Omega_{t,w}
\end{array}
\end{equation}
where $\mathbf{U_{w}}(t,x,y,z)=(U_{w},V_{w},W_{w})$,  $P_{w}(t,x,y,z)$, $\mathbf{F}$ stand for, respectively, 
the fluid velocity, the scalar pressure and the exterior strength of gravity.

A kinematic law at the free surface and  a no-leak condition $\mathbf{U_{w}}\cdot\mathbf{n}_{\textbf{wb}} = 0$ 
on the wet boundary complete the system 
($\mathbf{n}_{\textbf{wb}}$ being the outward unit normal vector, see \resim \ref{OyOz}). 

To model the effects of the air on the water, we write the  continuity of the normal stress at the boundary layer separating the two phases.
Thus, the water pressure $P_{w}$  may be written as: 
\begin{equation}\label{Pw}
P_{w} = g(h_{w}-z)\cos\theta + \frac{P_{a}}{\rho_0}
\end{equation}
where the term $g(h_{w}-z)\cos\theta$ is the hydrostatic pressure,  $\rho_0$ is the density of water at atmospheric pressure,
$P_{a}$ is the air pressure.
This is the second natural coupling between layers.

The Euler equations \eqref{EulerModelFluid} will be reduced to a shallow water like 
equations using standard arguments. \\
Let us briefly recall for the convenience of the reader the setup of the 
derivation (see \cite[Chapter 1]{TheseErsoy,BEG2012} for the detailed derivation).

The  aspect-ratio defined by 
$$\varepsilon = \frac{2R}{L}\ll 1 $$ is assumed small 
where $L$ and $R$ stands for the characteristic length and the  radius of the pipe. 
Following the classical thin-layer asymptotic analysis of Gerbeau and Perthame \cite{GP01} and  
\cite[Chapter 1]{TheseErsoy,BEG2012},  we introduce 
$\overline{U}$, $\overline{V}$, $\overline{W}$  respectively, the characteristic horizontal, normal and binormal 
speed with respect to the moving frame  such that 
$\epsilon = \frac{\overline{V}}{\overline{U}}= \frac{\overline{W}}{\overline{U}}, 
\, \overline{U} = \frac{L}{T},\,  P = \overline{U}^2$ for some characteristic time $T$. 

We also introduce the non dimensional variables
$\frac{U}{\overline U},\,  \varepsilon \frac{V}{ \overline U},\,  \varepsilon \frac{W}{\overline U},\, $
being respectively, the horizontal, normal and binormal non dimensional speed with respect to the moving frame with
$ \frac{x}{L},\,  \frac{y}{H},\,  \frac{z}{H},\,  \frac{p}{P},\, \frac{\theta}{\theta_0}$ 
respectively, the non dimensional $x$, $y$, $z$ coordinate,  the non dimensional pressure and the 
non dimensional angle.

Rescaling the Euler equations \eqref{EulerModelFluid} with these non dimensional variables, one can formally 
take $\varepsilon = 0$ to get the hydrostatic approximation.    
Then, we introduce the  conservative variables $A(t,x)$ and $Q(t,x)=A(t,x)u(t,x)$  
representing respectively the wet area and the discharge  defined as: 
$$A(t,x) = \int_{\Omega_{w}} dydz,\quad Q(t,x)=A(t,x)u(t,x)$$
where $u$ is the mean value of the speed 
 $$\dsp u(t,x)=\frac{1}{A(t,x)}\int_{\Omega_{w}} U_{w}(t,x,y,z)\,dydz.$$
As done by the authors in \cite{BEG2012}, taking the  averaged values  along the cross-section $\Omega_{w}(t,x)$ 
in the hydrostatic approximation of the incompressible Euler equations,
and approximating the averaged value of a product by the product of the averaged values, yields 
to the free surface model:
\begin{equation}\label{FSAirModel}
\left\{\begin{array}{lll}
\partial_t A +\partial_x Q &=& 0 \ ,\\
\dsp \partial_t Q+\partial_x\left(\frac{Q^2}{A}+A P_{a}/\rho_{0} +g I_1(x,A) \cos\theta
\right)&=&  -g A \partial_{x}Z \\
       & & +g I_2(x,A) \cos\theta\\[0.25cm]
       & & + \dsp P_{a}/\rho_{0}\, \partial_x A \ .
\end{array}\right.
\end{equation}
where $g$ is the gravity constant and $Z(x)$ is the elevation of the point $\omega(x)$. The terms $I_1(x,A)$ and $I_2(x,A)$ are defined by:  $$I_1(x,A) =
\int_{-R}^{h_{w}}(h_{w}-Z) \sigma(x,z) \,dz$$ and $$I_2(x,A) =
\int_{-R}^{h_{w}} (h_{w}-Z) \partial_x\sigma(x,z) \,dz\,.$$  
They represent respectively the  term of hydrostatic pressure and the pressure source term. In theses formulas $\sigma(x,z)$ is the width of the cross-section at position $x$ and at height $z$ (see \resim \ref{OyOz}).
\begin{rem}\label{remtheta2}
In view of Remark \ref{remtheta1}, whenever $\theta=\theta(x)$, an additional curvature source term of the form 
$\bar z (\cos\theta)'$ appears where $\bar z = h_w(A) - \frac{I_1(x,A)}{A}$ 
stands for the distance separating the free surface to the center of the mass of the fluid 
(see \cite[Chapter 1]{TheseErsoy,BEG2012}).
\end{rem}

\subsection{The air layer model}
The derivation of the air layer model is based on the Euler compressible equations  
\eqref{EulerModelGas} which writes in Cartesian coordinates $(x,y,z)$:
\begin{equation}\label{EulerModelGas}
\begin{array}{llll}
\partial_t \rho + \div (\rho \mathbf{U_{a}}) &=& 0, & \textrm{ on } \R\times\Omega_{t,a} \\
\partial_t (\rho \mathbf{U_{a}}) + \div (\rho \mathbf{U_{a}} \otimes \mathbf{U_{a}}) + \nabla P_{a} &=&0, & \textrm{ on }\R\times\Omega_{t,a}
\end{array}
\end{equation}
where $\mathbf{U_{a}}(t,x,y,z)=(U_{a},V_{a},W_{a})$ and $\rho(t,x,y,z)$ denotes the velocity  and  the density of the air whereas  $P_{a}(t,x,y,z)$ is the scalar pressure. Let us mention here that we neglect the gravitational effect on the air layer.

We assume  that the air layer  is isentropic, isothermal and follows a perfect gas law. Thus we have the following equation of state:
\begin{equation}\label{LaplacePressure}
\dsp P_{a}(\rho) = k_{a}\, \rho^{\gamma} \textrm{ with } \dsp k_{a} = \frac{p_{a}}{\rho_{a}^{\gamma}}
\end{equation}
for some reference pressure $p_{a}$ and density $\rho_{a}$. The adiabatic index $\gamma$ is set to $\dsp7/5$ for a diatomic gas and the constant $ k_{a}$
depends on the universal constant $R_{uc}$ and the temperature. 
\begin{rem}
The definition of the free surface pressure \eqref{Pw} and the air pressure 
\eqref{LaplacePressure}
ensures the continuity property of the normal stress for two immiscible perfect fluids at the interface.
\end{rem}
Taking advantage of the geometric configuration as in the previous subsection, i.e. the aspect-ratio 
$\varepsilon =2R/L$ being small, we proceed to a thin-layer asymptotic approximation following 
the derivation of the pressurized model (used in the context of the unsteady mixed flows in closed water pipes) 
\cite[Chapter 1]{TheseErsoy,BEG2012} which has the same structure as the 
Euler equations \eqref{EulerModelGas} used here. Indeed, only the index $\gamma$ changes from one model to this one.

The thin-layer asymptotic settings are essentially the same as used in the 
previous subsection. 
Thus, in the following development,  we  assume that the 
density of air and the pressure depends only on $t \mbox{ and } x$, i.e. we use this notation when 
no ambiguity occurs.

Let us now introduce the air area $\mathcal A$ by:
$$\mathcal A = \int_{\Omega_{a}} dydz \ ,$$
the averaged air velocity $v$ (we recall that $u$ is the averaged velocity of the fluid layer) by:
$$\dsp v(t,x)=\frac{1}{\mathcal A(t,x)}\int_{\Omega_{a}} U_{a}(t,x,y,z)\,dy dz \ ,$$
and the averaged density by:
$$\dsp \overline\rho(t,x)=\frac{1}{\mathcal A(t,x)}\int_{\Omega_{a}} \rho(t,x,y,z)\,dy dz \ .$$

The conservative variables are $M = \overline\rho/\rho_{0} \mathcal A$ (``pseudo area occupied by the air'') and $D = M v$ (``pseudo air discharge'').
Let us remark that we use rescaled mass and discharge instead of the real mass and discharge to be consistent with the free surface model \eqref{FSAirModel}.
As done before to obtain the water layer model, we take averaged values in the Euler equations \eqref{EulerModelGas} over sections $\Omega_{a}(t,x)$, 
and perform the same approximations on averaged values of a product to get the following model:

\begin{equation}\label{AirModel}
\left\{\begin{array}{lll}
\partial_t M + \partial_x D & = &
\dsp \int_{\partial \Omega_{a}} \rho/\rho_{0}\, \left(\partial_t \mathbf{m} +v\partial_x \mathbf{m} -   \mathbf{v}\right).  \mathbf{n}\, ds \ , \\
\partial_t D +
\partial_x\left(\dsp \frac{D^2}{M}+P_{a}(\overline{\rho})/\rho_{0}\,\mathcal{A}\right) &=&
\dsp P_{a}(\overline{\rho})/\rho_{0}\, \partial_x(\mathcal{A}) \\
 &  &+   \dsp \int_{\partial \Omega_{a}} \rho/\rho_{0}\, v \left(\partial_t \mathbf{m} +v\partial_x \mathbf{m} -   \mathbf{v}\right).  \mathbf{n}\, ds \ .
\end{array}\right.
\end{equation}
where $\mathbf{v} $ is  the velocity  in the $(\mathbf{j}, \mathbf{k})$-plane.
The boundary $\partial\Omega_{a}$ is divided into  the free surface boundary $\Gamma_{fs}$  and the  border of the pipe in contact  with  air $\Gamma_{c}$. For $m \in \Gamma_{c}$, $\dsp \mathbf{n} = \frac{\mathbf{m}}{|\mathbf{m}|}$ is the outward unit vector at the point $m$ in the $\Omega$-plane and $\mathbf{m} $ stands for the vector \textbf{${\omega m}$} while $\mathbf{n}$ denotes $-\mathbf{n}_{fs}$ on $\Gamma_{fs}$.
As the pipe is infinitely rigid, the non-penetration condition holds on $\Gamma_{c}$, namely: $\mathbf{U_{a}} \cdot\mathbf{n}_{ab}=0$ (see \resim \ref{OxOz}).

On the free surface, as the  boundary kinematic condition used for the water layer at the free surface holds, the integral appearing in System \eqref{AirModel} vanishes. Using the equation of state \eqref{LaplacePressure}, and 
using the air sound speed defined by:
\begin{equation}\label{airsound}
c_{a}^2 = \dsp \frac{\partial p } {\partial \rho} = k_{a}\gamma \dsp\left(\frac{\rho_0 M}{\mathcal {A}}\right)^{\gamma-1},
\end{equation}
we  finally  get the air-layer model:
\begin{equation}\label{AirlayerModel}
\left\{
\begin{array}{rcl}
\partial_t M + \partial_x D & = & 0 \ , \\
\partial_t D + \partial_x\left(\dsp \frac{D^2}{M}+ c_{a}^2 \frac{M}{\gamma} \right) &=& \dsp c_{a}^2 \frac{M}{\gamma}  \, \frac{\partial_x(\mathcal A)}{\mathcal A} \ .
\end{array}\right. 
 \end{equation}
\subsection{The two-layer model}\label{SubsectionBilayerModel}
The two-layer averaged flow model is  then simply obtain by the apposition of the models \eqref{FSAirModel} and \eqref{AirlayerModel}, using the coupling  \eqref{AirWaterRelation} 
which also writes  $\mathcal A+A = S$ where $S=S(x)$ denotes the pipe section and 
using the definition of the air pressure \eqref{LaplacePressure} and  the air sound speed \eqref{airsound},
we get:
\begin{equation}\label{BilayerModel}
\left\{\begin{array}{rcl}
\partial_t M + \partial_x D & = & 0 \ , \\
\partial_t D + \partial_x\left(\dsp \frac{D^2}{M}+  c_{a}^2 \frac{M}{\gamma}\right) &=& \dsp  c_{a}^2 \frac{M}{\gamma}\,\frac{\partial_x(S-A)}{(S-A)} \ , \\
 & & \\
\partial_t A +\partial_x Q &=& 0 \ ,  \\
\dsp \partial_t Q+\partial_x\left(\frac{Q^2}{A}+g I_1(x,A) \cos\theta
+ A \frac{c_{a}^2 \, M}{\gamma \, (S-A)}  \right)&=&  -g A \partial_{x} Z+g I_2(x,A) \cos\theta\\
       & & + \dsp \frac{c_{a}^2 \, M}{\gamma \, (S-A)} \, \partial_x A \ .
\end{array}\right. 
\end{equation}
\subsection{Study of the eigenvalues of the  two-layer model}
As for almost every two-layer shallow water models \cite{CMP01,S07,BM08,AK09} 
as well as two phase fluid flows \cite{I75,S76,SW84,HI03,Sainsaulieu93}, our model is a 
\emph{conditionally hyperbolic} system as we will see thereafter. By definition, it means that with respect to some physical parameters,
for instance, when the relative speed between the two layers is small the system remains hyperbolic as already pointed out 
by several authors (see for instance \cite{Sainsaulieu93,BM08,AK09}). In fact, 
this property is not only restricted to small relative speed but holds also for large one as it was 
emphasized by Ovsjannikov \cite{Ovsjannikov79} and recalled later on by Barros and Choi \cite{BC08} and the reference
 therein.

In the present case, the coupling between the layers, due to the hydrostatic pressure \eqref{Pw} 
and the barotropic one \eqref{LaplacePressure}, provides  a full access to the 
the computation of the eigenvalues, which is required for almost every numerical schemes based on Riemann solver. 

Noting  $\mathbf{W} = (M,D,A,Q)^t$ the unknown vector,  the system may be written in the quasilinear form:
\begin{equation*}
\partial_t \mathbf{W}+\mathcal{D}(x,\mathbf{W}) \partial_X \mathbf{W} =0
\end{equation*}
with the  convection matrix:
$$ \mathcal{D}= \left(
\begin{array}{cccc}
0 & 1 & 0 & 0 \\
c_{a}^2-v^2 & 2v & \dsp \frac{M}{S-A}c_{a}^2 & 0 \\
0 & 0 & 0 & 1 \\
\dsp\frac{A}{(S-A)}c_{a}^2 & 0 & c_{w}^2+\dsp\frac{AM}{(S-A)^2}c_{a}^2-u^2 & 2u
\end{array}\right)$$
contains flux gradient terms for which 
the quantities  $c_{a} $ stands for the air sound speed defined by \eqref{airsound} and  
$c_{w} =  \dsp\sqrt{g \frac{A}{T(A)}\cos\theta}$  the water sound speed where $T(A)$ is 
the width of the free surface at height $h_{w}(A)$ (see \resim \ref{OyOz}). 
\begin{rem}
Let us emphasize that $c_w$ stands for the classical water sound speed when the air layer is not taken into account. The new quantity 
$c_m^{2}:=c_{w}^2+ \frac{AM}{(S-A)^2}c_{a}^2$ represents the square of the water sound speed which takes into account the air effect. 
\end{rem}

The characteristic polynomial of degree $4$ of the convection matrix  $\mathcal{D}$ of the two-layer system is defined by:
\begin{equation}\label{QuarticTwoLayer}
P(\lambda)  = 
\dsp\left(\lambda^2-2v\lambda -\left(c_{a}^2-v^2\right)\right)
\left(\lambda^2-2u\lambda-\left(c_{w}^2+\frac{AM}{(S-A)^2}c_{a}^2-u^2\right)\right) \\
- \dsp \frac{{AM}}{(S-A)^2} c_{a}^4
\end{equation}
for which no simple expression for the eigenvalues can be computed but  
one can provide, as it was already done in the context of the two-layer shallow water model, 
the first order expansion with respect to the relative speed 
between layers (see for instance \cite{schijf1953theoretical,A05}) 
and derive an  hyperbolicity condition which may also linked to the Kelvin-Helmholtz instability. 

In what follows, we define the following non dimensional numbers: 
\begin{equation}\label{VarNonDim}
F = \frac{v-u}{c_m},\quad \sqrt{H} = \frac{c_a}{c_m},\quad 
c_m=\sqrt{c_w^2+s c_a^2} \textrm{ with } s= \frac{AM}{(S-A)^2}\geqslant 0\, .
\end{equation}
$F$ stands for the relative speed between the two layers, $H$ represents the ratio between the sound speed of each layer and $s = \frac{\rho}{\rho_0}\frac{A}{S-A}$
represents the compressibility ratio between the air layer and the water layer.

In Section \ref{KineticFramework}, based on the kinetic formulation of the two-layer model \eqref{BilayerModel}, we introduce a 
new numerical scheme which does not require the computation of the eigenvalues and therefore we will focus  on its stability and 
sensitivity when the system is non hyperbolic. This motivates the following study. 

\begin{rem}
Let us firstly remark that although, for the construction of the kinetic scheme we are not interested in the the computation of these eigenvalues, 
one can provide easily the expression  of them when the relative speed is equal to zero i.e. $F = 0$,: 
$$
\dsp u \pm
\frac{1}{2}
\sqrt{  2 c_a^2+ 2 c_m^2 \pm \sqrt{ \left(c_a^2-c_m^2\right)^2 +4 s c_a^4 } }
$$
for which we immediately see that, we may have complex eigenvalues as  observed on \resim\ref{NH1000AFDeZoom} and \ref{NH1000AF}. 
When  $H \geq 1/4s$ eigenvalues are real, otherwise  two eigenvalues are real and two are complex.  
Moreover, in comparison with the two-layer shallow water for instance, 
the condition ``relatively small enough'' is not sufficient to state that the present system is hyperbolic. 
This is due to the fact that $s$ is not bounded from above (as we will see). 
\end{rem}
Let us recall the so-called \emph{root location criteria} for a polynomial of degree 4 (see also \cite{Fuller81}).
\begin{thm}\label{thmquartic}
Let $P$ be a polynomial of degree 4 defined by:
\begin{equation}\label{quartic}
P(x) = \sum_{k=0}^{4} a_{k} x^{4-k}  \textrm{ for } (a_k)_{k} \in \R \textrm{ and } a_0>0. 
\end{equation}
All the root of Equation  \eqref{quartic} are real if and only if one of the following conditions holds:
\begin{enumerate}
 \item\label{cond1} $\Delta_3> 0$, $\Delta_5> 0$ and $\Delta_7 \geqslant 0$, 
 \item\label{cond2}$\Delta_3\geqslant 0$, $\Delta_5=0$ and $\Delta_7 = 0$.
\end{enumerate}
where $\Delta_3$,  $\Delta_5$,  $\Delta_7 $ are the determinant of the matrices:
$$
\Delta_3 = 
\left(
\begin{array}{ccc}
a_0    &  a_1   & a_2\\
0      &  4 a_0 & 3 a_1 \\ 
4a_0   &  3a_1  & 2 a_2 
\end{array}
\right), \quad
\Delta_5 = 
\left(
\begin{array}{ccccc}
a_0  & a_1    &  a_2   & a_3     &  a_4 \\
0    & a_0    &  a_1   & a_2     &  a_3 \\
0    & 0      &  4 a_0 & 3 a_1   &  2a_2 \\ 
0    & 4a_0   &  3a_1  & 2 a_2   &  a_3\\
4a_0 & 3a_1   &  2a_2  & a_3     &  0 \\
\end{array}
\right), 
$$
and 
$$
\Delta_7 = 
\left(
\begin{array}{ccccccc}
a_0  & a_1  &  a_2   & a_3    &  a_4    & 0      & 0 \\
0    & a_0  & a_1    &  a_2   & a_3     &  a_4   & 0  \\
0    & 0    & a_0    &  a_1   & a_2     &  a_3   & a_4  \\
0    & 0    & 0      &  4 a_0 & 3 a_1   &  2a_2  & a_3 \\ 
0    & 0    & 4a_0   &  3a_1  & 2 a_2   &  a_3   & 0 \\
0    & 4a_0 & 3a_1   &  2a_2  & a_3     &  0     & 0 \\
4a_0 & 3a_1 & 2a_2   & a_3    &  0    & 0      & 0 \\
\end{array}
\right).
$$
\end{thm}
Let us now apply this result to the characteristic equation \eqref{QuarticTwoLayer}. To this end, we closely follow the study of Barros and Choi \cite{BC08}. 
We  thus introduce the non-dimensional variable $x = \lambda/c_m$ and as mentioned by Barros and Choi, without loss of generality, 
we may assume that $u/c_m = 1$ (by choosing a moving reference frame such that this condition is met). Then 
Equation \eqref{QuarticTwoLayer} becomes:
$$P(x) = x^4 - 2(2+F)x ^3 + \left((1+F)(5+F)-H\right)x^2+2\left(H-(1+F)^2\right)x -s H^2.$$ 

Let us first remark that $\Delta_3 =4 (F^2+2(1+H))$ is always strictly positive.

Secondly, we have:
$$\Delta_5=8
\Big(
(1+H)F^4-2\left(H^2(1-s)+1-6H\right)F^2 +(1+H)\left((H-1)^2+4 s H^2\right)
\Big).$$
Thus $\Delta_5$ is a polynomial of degree 2 in the variable $y=F^{2}$.
Its  discriminant $R(H,s)$ is equal to:
$$R(H,s)= 256H \left( (s^2-6s)H^3+(4s-12)H^2+(40-6s)H-12 \right) $$ 
which is a polynomial of degree 4 in the variable $H$. Since $8(1+H)$ (which is the dominant coefficient of $\Delta_{5}$) is positive, we are  firstly interested
in the set where the discrimant $R(H,s)$ is negative.

Since $s\geqslant 0$,  $R(H,s)$  admits four real roots :
$$l_1(s) = \frac{-1+\sqrt{1+6s}}{s},\,\,\,   l_2(s) = -\frac{1+\sqrt{1+6s}}{s},\,\,\,   l_3(s) = -\frac{2}{s-6}, \textrm{ and } l_4(s) =0.$$
Let  $\mathcal{D}$ bet the set defined by:
$$\mathcal{D} = 
\left\{
(H,s)\in \R^2_+; 
\begin{array}{l}
0<s<4 \textrm{ and } l_3(s)<H<l_1(s),\\
\textrm{ or } 4<s<6 \textrm{ and } l_1(s)<H<l_3(s),\\
\textrm{ or } s>6 \textrm{ and } H>l_1(s) 
\end{array}
\right\}$$
If  $(H,s)\in\mathcal{D}$, then $R(H,s)\leqslant 0$ and  $\Delta_5>0$. 
On the contrary, if  $(H,s)\notin\mathcal{D}$ then $R(H,s)\geqslant 0$ and $\Delta_5$  admits two real eigenvalues which are both negative:
$$y_{\pm} = \frac{2\left(H^2(1-s)+1-6H\right) \pm \sqrt{R(H,s)}}{(1+H)}<0 \textrm{ for } (H,s)\notin \mathcal{D}  .$$
By definition, we have $y = F^{2}$ and thus for every $(H,s) \in \R^{2} \,,\, \Delta_5 > 0$. Therefore we must only examine the sign of $\Delta_{7}$.

$\Delta_7$ is a polynomial of degree 4 in the variable $y=F^2$ and we write it as $\Delta_7(y) = 16 H Q(y)$ where 
$$
\begin{array}{lll}
Q(y) &=& y^4 + \left(sH^2+(s-4)H-4\right)y^3 +\left( (s^2-3s)H^3 + (6-26s)H^2 +(4-3s)H+6\right)y^2\\
     & & +\left( (3s-20s^2)H^4 + (13s-20s^2-4)H^3 +(13s+4)H^2+ (4s+3)H-4\right)y\\ 
     & & -(16s^3+8s^2+s)H^5 + (32s^2+12s+1)H^4\\
     & & -(4+22s+8s^2)H^3+(12s+6)H^2-(4+s)H+1\, .
\end{array}
$$
We have to study the sign of $\Delta_{7}(y)$ only for $y \geq 0$.  Using a symbolic computation software 
 such as  Maxima$^\copyright$, we obtain:
 $$ \Delta_3(Q)>0,\,\Delta_5(Q)\leqslant 0  \mbox{ and }\Delta_7(Q)>0 \ .$$ 
Applying again Theorem \ref{thmquartic} (see also Fuller \cite[Theorem 2 and 4]{Fuller81}) for the polynomial $Q$,
it exists $y_{min}(H,s) \in \R$ which may be positive or negative and $y_{max}(H,s) > 0$ (which are the  two real roots of $Q$) such that:
$$\forall y \in (-\infty,y_{min}]\cup [y_{max},\infty) \,,\, Q(y) \geq 0 \, .$$ 

\begin{figure}[H]
\subfigure[$y_{min} <0$\label{quarticcond}]{\includegraphics[height = 4.0001cm]{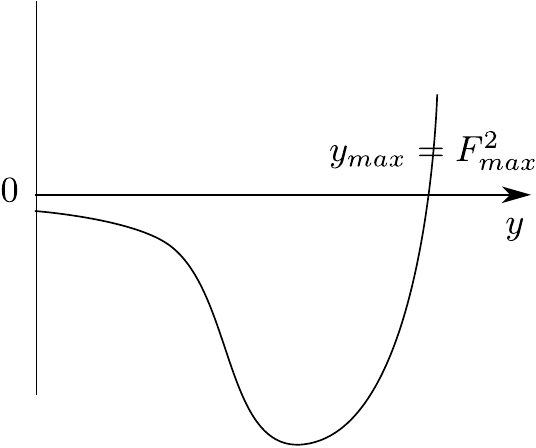}}
\subfigure[$y_{min}>0$\label{quarticIncond}]{\includegraphics[height = 4.0001cm]{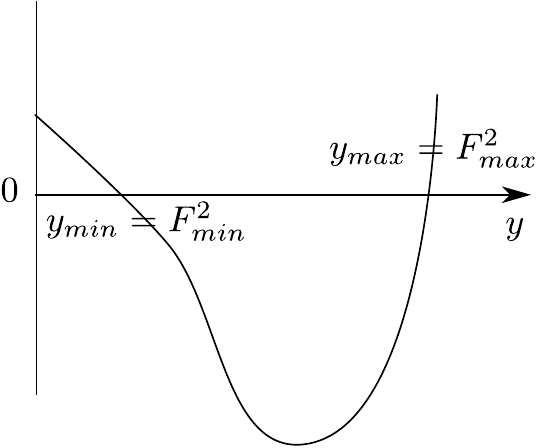}}
\caption{Behavior of the polynomial $Q(y)$\label{quarticond}} 
\end{figure}

We have illustrated this property in  \resim\ref{quarticond}.

We can now summarized the study of the hyperbolicity of the  two-layer model by the following result:
\begin{thm}\label{hyperbolicity}
Given a couple $(s,H)$, the two-layer system is hyperbolic whenever $y= F^{2}$ satisfies one of the following conditions:
\begin{enumerate}
 \item $y_{min}(H,s)\leqslant 0 \mbox{ and  } y \geqslant y_{max}$
 \item $y_{min}(H,s)>0\mbox{ and  } y \in [0,y_{min}] \cup [y_{max},\infty)$
\end{enumerate}
\end{thm}
\begin{rem} Let us conclude this study about the hyperbolicity of the two layer model by the following remarks.
\begin{enumerate}
\item  The criterion $y= F^2\geqslant F^2_{max}:=y_{max}$ ensures the hyperbolicity for large relative speed while when $y_{min}>0$
the condition $0 \leqslant F^2 \leqslant y_{min}=F^2_{min}$  corresponds to small relative speed. 
\item The term $s=\frac{\rho}{\rho_0}\frac{A}{S-A}$ is the equivalent ratio $\frac{\rho_2}{\rho_1}<1$ appearing in the two-layer shallow water equations, see \cite{BC08}.
Here $s$ is not  bounded from above thus  $y_{min}$ is not necessary positive.
\item Since $s=s(\rho,A)$ and $H =H(\rho,A)$, the preceding result can be formulated in terms of the following three parameters $F$, $\rho$ and $A$.
For different values of these parameters, we may draw the region where one of the preceding conditions is fulfilled.\\
In \resim\ref{NHCondHyperbolic1} and \ref{NHCondHyperbolic2} the black gray represents the non hyperbolic region.  \\
In \resim\ref{NH1000AFDeZoom} and \ref{NH1000AF}, we represent the hyperbolic region part for small 
relative speed between the two layers (that is for large value of $\rho$) while in \resim\ref{NH10AF} we represent the hyperbolic region part for large
relative speed between the two layers (that is for small value of $\rho$). We see in \resim\ref{NH10AF} 
that the thickness of the hyperbolic region decreases with $\rho$. Finally we also represent on \resim\ref{NH} the hyperbolic region versus the three parameters
$F$, $\rho$ and $A$.
\begin{figure}[H]
\subfigure[$\rho=1000$,  $F\: x$-axis and $A\: y$-axis\label{NH1000AFDeZoom}]
{
\includegraphics[width= 6.6cm]{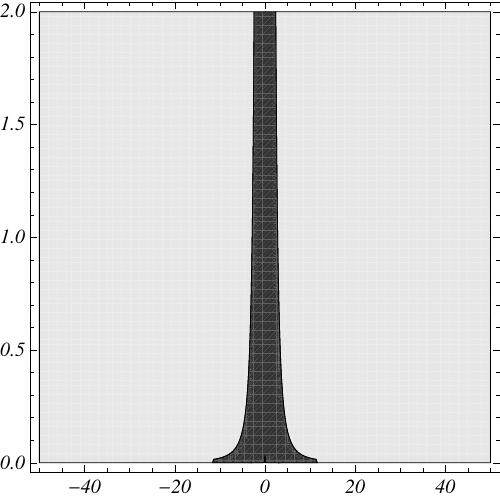}
} 
\subfigure[Zoom on :$\rho=1000$,  $F\: x$-axis and $A\: y$-axis\label{NH1000AF}]
{
\includegraphics[width = 6.6cm]{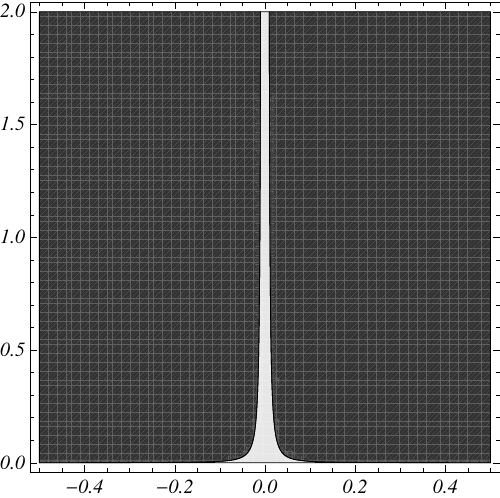}
} 
\caption{Small and large relative speed: hyperbolic region \label{NHCondHyperbolic1}}
\end{figure}

\begin{figure}[H]
\subfigure[$\rho=10$,  $F\: x$-axis and $A\: y$-axis\label{NH10AF}]
{
\includegraphics[width = 6.6cm]{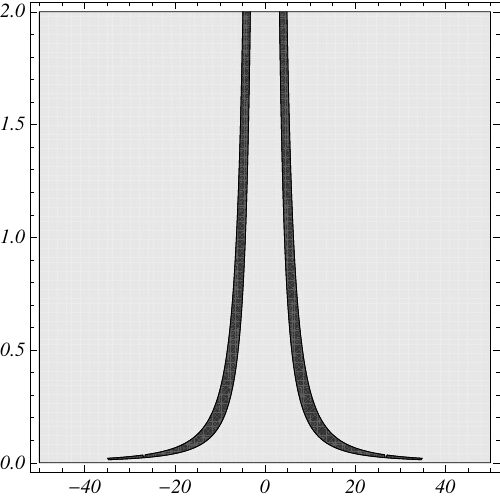}
} 
\subfigure[3d hyperbolic region for :$0\leqslant  \rho \leqslant 1000$ ($x$-axis),  $-2\leqslant  F \leqslant 2$  ($y$-axis), 
$0\leqslant  A \leqslant 2$  ($z$-axis)\label{NH}]
{
\includegraphics[width = 6.6cm]{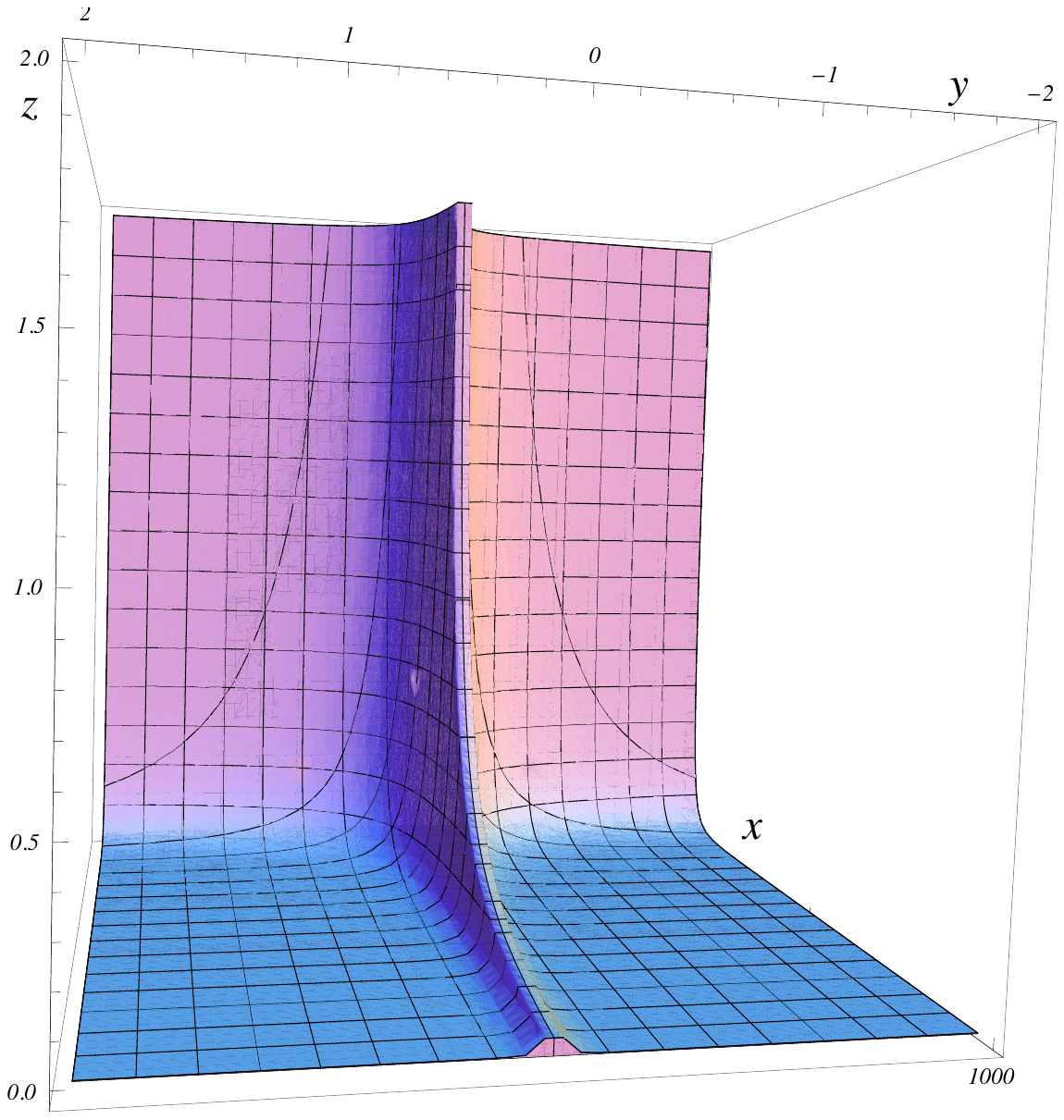}
} 
\caption{Strip hyperbolic range and the hyperbolic set\label{NHCondHyperbolic2}}
\end{figure}
\item We conclude the preceding numerical study of the eigenvalues by noticing that for small and large relative speed, 
the two layer model remains hyperbolic.  
\end{enumerate}
\end{rem}

Finally, even if  System \eqref{BilayerModel} is conditionally hyperbolic, 
the existence of a convex entropy function is not a contradiction. 
Consequently, admissible weak solutions should satisfy the corresponding below-mentioned entropy inequality 
(see for instance \cite{BM08,AK09}). 

The results are summarized in the following theorem:
\begin{thm}\label{BilayerModelThm}~
\begin{enumerate}
 \item For smooth solutions of (\ref{BilayerModel}), the velocities $u$ and $v$ satisfy
    \begin{eqnarray*}
          \partial_t v+\partial_x\left(\dsp \frac{v^2}{2}+\frac{c_{a}^2}{\gamma-1}\right)&=& 0  \ , \\
          \partial_t u +\partial_x\left(\dsp \frac{u^2}{2}+\dsp g h_w(A)\cos\theta + gZ+ \frac{ c_{a}^2 \, M}{\gamma(S-A)}\right)&=&0  \ .
    \end{eqnarray*}
\item System \eqref{BilayerModel} admits a  mathematical total energy $$\mathcal{E} = E_{a} + E_{w}$$
with $$E_{a} = \frac{Mv^2}{2}+\frac{c_{a}^2 \, M}{\gamma(\gamma-1)}$$
and
 $$E_{w} = \frac{Au^2}{2}+g A (h_w-I_1(x,A)/A)\cos\theta + g A Z$$
which satisfies, for smooth solutions, the entropy  equality $$\partial_t  \mathcal{E} +\partial_x \mathcal{H} = 0$$ where the flux $\mathcal{H}$ 
is the total entropy flux defined by 
$$\dsp \mathcal{H} = \mathcal{H}_{a} +\mathcal{H}_{w}$$
where the air entropy flux reads:
$$\mathcal{H}_{a} =  \left(E_{a}+\frac{c_{a}^2 \, M}{\gamma}\right)v $$
and the water entropy flux reads: 
$$\mathcal{H}_{w} =  \left(E_{w} + gI_1(x,A)\cos\theta+A\frac{c_{a}^2 \, M}{(S-A)}\right)u \ .$$
\item The energies $E_{a}$ and $E_{w}$ satisfy the following entropy flux equalities:
\begin{equation}\label{AirEnergy}
\partial_t E_{a} + \partial_x \mathcal{H}_{a} = \frac{c_{a}^2\,M}{\gamma(S-A)} \, \partial_t A
\end{equation}
and
\begin{equation}\label{WaterEnergy}
\partial_t E_{w} + \partial_x \mathcal{H}_{w} =- \frac{c_{a}^2\, M}{\gamma(S-A)} \, \partial_t A \ .
\end{equation}
\end{enumerate}
\end{thm}
\begin{proof}[Proof of Theorem \ref{BilayerModelThm}]
The proof of these statements relies only on algebraic combinations of the two equations forming System \eqref{BilayerModel}.
\end{proof}
\begin{rem}
Let us remark that Equations \eqref{AirEnergy}-\eqref{WaterEnergy} can be interpreted in such a way:
the system is physically closed i.e. the energy dissipated (or gained) by the water layer is forwarded to (or lost by)  the air layer. 
This property ensures the total energy equality $$\partial_t \mathcal{E} +\partial_x \mathcal{H} = 0\ .$$
\end{rem}
\section{The kinetic interpretation of the two-layer air-water model}\label{Kinetic}
As already mentioned before, the coupling between the layers, due to the coupling of an hydrostatic pressure and a barotropic 
pressure law, does not provide an explicit access to the eigenstructure of the system, 
which is required for numerical schemes based on the solution of the Riemann problem. It forms a system which is conditionally hyperbolic 
(see Theorem \ref{hyperbolicity}). Nevertheless, one can provide the first order expansions with respect to the relative speed as approximate expressions for the eigenvalues 
(see for instance \cite{schijf1953theoretical,A05,AK09}).  
Let us mention that Abgrall and Karni \cite{AK09} provides a way to decouple 
the system and to have an access to the eigenvalues by using a relaxation approach. 

In this work, the computation of the eigenvalues is  out of interest since we will focus on the 
kinetic interpretation of the two-layer model. This is a useful way to compute the numerical solution even if the system loses its hyperbolicity since 
the computation of the eigenvalues are not required. 
Moreover following \cite{BEG11_1,BEG11_2},  natural properties are obtained such as the apparition or vanishing of vacuum in the air layer 
or drying and flooding flows for the water layer. All of these properties will be useful in the future development for more physical models. 
Although at the numerical level, the effect of the variable cross section is not taken into account, this framework offers an easy way to 
upwind the source terms and to get well-balanced numerical schemes. 

The particular treatment of the boundary conditions will be rapidly exposed using the ``decoupled system'' approach.
\subsection{The mathematical kinetic formulation}\label{KineticFramework}
Let $\chi:\R\to\R$ be a given real function satisfying the following  properties:
\begin{equation}\label{propchi}
\chi(\omega)=\chi(-\omega) \geqslant 0\;,\;
\int_{\R} \chi(\omega) d\omega =1,
\int_{\R} \omega^2 \chi(\omega) d\omega=1 .
\end{equation} 
It permits to  define  the density of particles, by a so-called \emph{Gibbs equilibrium}, 
\begin{equation}\label{Gibbs}
\mathcal{M}_{\alpha}(t,x,\xi) =
\frac{A_{\alpha}(t,x)}{b_{\alpha}(t,x)} \chi\left(\frac{\xi-u_{\alpha}(t,x)}{b_{\alpha}(t,x)}\right)
\end{equation}
where $b_{\alpha}(t,x) = b(x,A(t,x),M(t,x))$ with 
\begin{equation*}
b_{\alpha}^2(x,A,M) =
\left\{
\begin{array}{lll}
\dsp \frac{c_{a}^2}{ \gamma M} & \textrm{ if } & \alpha = a \ , \\
~\\
\dsp g\frac{I_1(x,A)}{A}\cos\theta + \frac{c_{a}^2  \, M}{\gamma \, (S-A)} & \textrm{ if } & \alpha = w \ ,
\end{array}\right. 
\end{equation*}
and
$$
\begin{array}{l}
A_{\alpha} =
\left\{
\begin{array}{lll}
A & \textrm{ if } & \alpha = w\\
M & \textrm{ if } & \alpha = a
\end{array}
\right.  , \
Q_{\alpha} =
\left\{
\begin{array}{lll}
Q & \textrm{ if } & \alpha = w\\
D & \textrm{ if } & \alpha = a
\end{array}\right. \mbox{ and } \
u_{\alpha} =
\left\{
\begin{array}{lll}
u & \textrm{ if } & \alpha = w\\
v & \textrm{ if } & \alpha = a 
\end{array}
\right. \ .
\end{array}
$$

The Gibbs equilibrium $\mathcal{M}_{\alpha}$ is related to the two-layer  model \eqref{BilayerModel} by the classical 
 \emph{macro-micro}scopic kinetic relations:
\begin{equation}\label{macroA}
A_{\alpha}  = \dsp\int_{\R} \mathcal{M}_{\alpha}(t,x,\xi)\,d\xi\,,
\end{equation}
\begin{equation}\label{macroQ}
Q_{\alpha}  = \dsp\int_{\R} \xi\mathcal{M}_{\alpha}(t,x,\xi)\,d\xi\,,
\end{equation}
\begin{equation}\label{macroFlux}
\dsp\frac{Q_{\alpha}^2}{A_{\alpha} }+A_{\alpha}\,b_{\alpha}(x,A,M)^2  = \dsp\int_{\R} \xi^2 \mathcal{M}_{\alpha}(t,x,\xi)\,d\xi\,.
\end{equation}
From the relations \eqref{macroA}--\eqref{macroFlux}, the nonlinear two-layer model can be viewed as two single 
linear equations, one for each layer, involving the nonlinear quantity 
$\mathcal{M}_{\alpha}$:
\begin{thm}[Kinetic Formulation of the \textbf{PFS} model]\label{ThmKineticFormulationPFS}
Assuming a constant angle $\theta$, 
$(A,Q,M,D)$ is a strong solution of System \eqref{BilayerModel} if and only if $(\mathcal{M}_{w,}\mathcal{M}_{a)}$ satisfies 
the system of coupled kinetic transport equations:
\begin{equation}\label{KineticFormulationPFS}
\left\{
\begin{array}{l}
\partial_t \mathcal{M}_{w}+\xi \cdot \partial_x\mathcal{M}_{w} - 
\left(g \partial_{x}Z  - g\dsp\frac{I_2(x,A)}{A} \cos\theta\,  -  \frac{M}{\gamma A (S-A)}c_{a}^2  \partial_{x}A\right)\,\partial_\xi \mathcal{M}_{w} = 
\mathcal{K}_{w}(t,x,\xi) \\
\partial_t \mathcal{M}_{a}+\xi \cdot \partial_x\mathcal{M}_{a} -  \left(- \, \frac{c_{a}^2}{\gamma(S-A)}\,\partial_{x}(S-A)\right)\partial_\xi \mathcal{M}_{a} = 
\mathcal{K}_{a}(t,x,\xi)
\end{array}
\right.
\end{equation}
for some collision term $\mathcal{K}_{w}(t,x,\xi)$ and $\mathcal{K}_{a}(t,x,\xi)$ which satisfies for $(t,x)$ a.e. 
$$\dsp  \int_{\R} \vecdeux{1}{\xi} \mathcal{K}_{w}(t,x,\xi)\,d\xi = 0 \;,\;  \int_{\R} \vecdeux{1}{\xi} \mathcal{K}_{a}(t,x,\xi)\,d\xi = 0.$$
\end{thm}
\begin{proof}
The proof relies on very obvious computations since $\mathcal{M}_{\alpha}$ verifies the macro-microscopic kinetic relations \eqref{macroA}, \eqref{macroQ}, \eqref{macroFlux}.
\end{proof}
\begin{rem}~
\begin{itemize} 
\item In order to write the system \eqref{KineticFormulationPFS} in a compact form, we use the indexes $w$ for the water and $a$ for the air and we
obtain:
\begin{equation*}
\partial_t \mathcal{M}_{\alpha}+\xi \cdot \partial_x\mathcal{M}_{\alpha} - \phi_{\alpha} \,\partial_\xi \mathcal{M}_{\alpha} = \mathcal{K}_{\alpha}(t,x,\xi)
\end{equation*}
The source terms $\phi_{\alpha}(x,\WW)$ are defined as:
\begin{equation}\label{PFSSourceTermPhi}
\phi_{\alpha}(x,\WW) = \BB_{\alpha}(x,\WW)\; ^{t}\big(\partial_x \WW_{\alpha}\big)
\end{equation}
with 
\begin{equation}\label{WW}
\WW_{\alpha} = 
\left\{
\begin{array}{lll}
\left(Z,S,A\right)    & \textrm{ if } & \alpha = w\ ,\\
\mathcal{A}:=S-A      & \textrm{ if } & \alpha = a\ ,
\end{array}
\right.
\end{equation}
and
$
\BB_{\alpha} =
\left\{
\begin{array}{ll}
\dsp \left(g \,,\, - g\dsp\frac{I_2(x,A)}{A \partial_x S} \cos\theta\,,\, -  \frac{M}{\gamma A (S-A)}c_{a}^2\right)  & \textrm{ if } \alpha = w,\\
\dsp - \, \frac{c_{a}^2}{\gamma(S-A)}  & \textrm{ if }  \alpha  = a \ .\\
\end{array}
\right.
$
\item The kinetic formulation presented in Theorem \ref{ThmKineticFormulationPFS} is a (non physical) microscopic
description of the two layer model \eqref{BilayerModel}. 
\item The results of Theorem \eqref{ThmKineticFormulationPFS} are obtained assuming a constant angle $\theta$. 
The generalization of these results with varying $\theta=\theta(x)$  can be obtained by adding the term of curvature mentioned in Remark \ref{remtheta2}
as done  in \cite{BEG09_2, BEG11_1, BEG11_2}, for a single fluid flow.
It can be easily adapted to this case.
\end{itemize}
\end{rem}

\section{The two-layer kinetic scheme}\label{KineticScheme}
In this section, following the works of \cite{BEG11_1,BEG11_2,PS01}, we will construct a two-layer finite volume kinetic scheme that 
is area conservative, that preserves the posivity of the area and that computes ``naturally'' flooding, drying and vacuum zones, even when the system loses its hyperbolicity.
Although, we are not interested in the treatment of the source term, the full kinetic scheme will be provided. 

Under this assumption and based on the kinetic formulation (see Theorem \ref{ThmKineticFormulationPFS}), we construct easily the two-layer  Finite
Volume scheme where the conservative quantities are cell-centered and  source terms are included into the numerical fluxes by a standard
kinetic scheme with reflections \cite{PS01}. 
\subsection{The general formulation}
\begin{figure}[H]
 \begin{center}
 \includegraphics[height = 2cm]{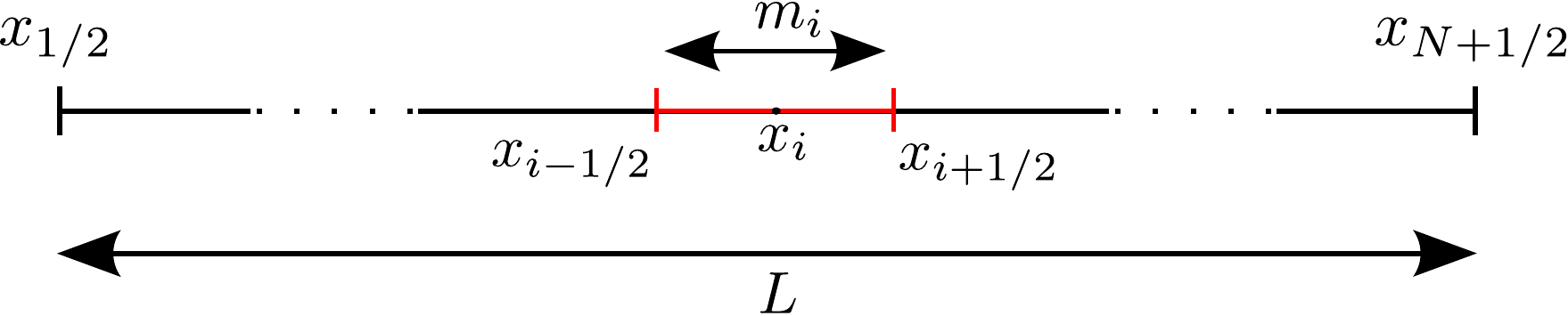}
  \caption{The space discretisation.}
  \label{pipediscret}
 \end{center}
\end{figure}

Let $N\in\N^{*}$, and let us consider the following mesh on $[0,L]$ of a pipe of length $L$. Cells are denoted for every
$i\in [0,N+1]$, by $m_i =(x_{i-1/2},x_{i+1/2})$, with  $x_i=\frac{x_{i-1/2}+x_{i+1/2}}{2}$ and 
$h_{i}=x_{i+1/2}-x_{i-1/2} $ the space step. The ``fictitious'' cells $m_0$ and $m_{N+1}$ denote the boundary cells and 
the mesh  interfaces located at $x_{1/2}$ and $x_{N+1/2}$ are respectively the upstream and the downstream ends of the pipe (see \resim\ref{pipediscret}).
We also consider a time discretization $t^n$ defined by $t^{n+1}=t^n+\Delta t^n$ with $\Delta t^n$ the time step. 

We denote $\UU_{\alpha,i}^n=(A_{\alpha,i}^n,Q_{\alpha,i}^n)$, $\dsp u_{\alpha,i}^n = \frac{Q_{\alpha,i}^n}{A_{\alpha,i}^n}$, $\mathcal{M}_{\alpha,i}^n$ 
the cell-centered approximation of $\UU_{\alpha} = (A_{\alpha},Q_{\alpha})$, $u_{\alpha}$ and $\mathcal{M}_{\alpha}$ on the cell $m_i$ at time $t^n$ where 
we use the notation subscript $\alpha$ introduced previously. 
We denote by $\UU_{\alpha,0}^n=(A_{\alpha,0}^n, Q_{\alpha,0}^n)$ the upstream  and $\UU_{\alpha,N+1}^n=(A_{\alpha,N+1}^n, Q_{\alpha,N+1}^n)$ 
the downstream state vectors.

The piecewise constant representation of $\dsp\WW_{\alpha}$ defined by \eqref{WW} is given by,   
$\WW_{\alpha}(t,x) = \dsp \WW_{\alpha,i}(t) \mathds{1}_{m_i}(x)$ 
where $\WW_{\alpha,i}(t)$ is defined as $\WW_{\alpha,i}(t) = \dsp \frac{1}{\Delta x}\int_{m_i}\WW_{\alpha}(t,x)\,dx$ for instance.

Denoting by $\WW_{\alpha,i}$ and $\WW_{\alpha,i+1}$ the left hand side and the right hand side values of $\WW$ at  the cell interface $x_{i+1/2}$, 
and using the ``straight lines'' paths (see  \cite{DLM95}) 
$$\Psi(s,\WW_{\alpha,i},\WW_{\alpha,i+1}) = s\WW_{\alpha,i+1}+(1-s)\WW_{\alpha,i},\,s\in[0,1],$$ 
we define the non-conservative product $\phi_{\alpha}(t,x_{i+1/2})$ by writing:
\begin{equation}\label{DefinitionOfNonConservativeProductPhi}
\phi_{\alpha}(t,x_{i+1/2})  = \jump{\WW_{\alpha}}(t)\cdot\int_0^1 \BB\left(t,\Psi(s,\WW_{\alpha,i}(t),\WW_{\alpha,i+1}(t))\right)ds
\end{equation}
where  $\jump{\WW_{\alpha}}(t):=\WW_{\alpha,i+1}(t)-\WW_{\alpha,i}(t)$, is the jump of $\WW_{\alpha}(t)$ across the discontinuity localized 
at $x=x_{i+1/2}$. As the first component of $\BB_w$ is $g$, we recover  the classical interfacial upwinding for the  conservative term $Z$ 
as initially introduced in \cite{PS01}.

Neglecting the collision kernel as in \cite{PS01,BEG11_1,BEG11_2} and using the fact that 
$\phi_{\alpha} = 0$ on the cell $m_i$ (since $\WW_{\alpha}$ is constant on $m_{i}$), the kinetic transport equation
\eqref{KineticFormulationPFS} simply reads:
\begin{equation}\label{eqcin}
\frac{\partial}{\partial t} {\mathcal M}_{\alpha} +\xi \cdot \frac{\partial}{\partial x} {\mathcal M}_{\alpha}
=0 \quad \mbox{ for } x \in m_i \ .
\end{equation}
This equation is a linear transport equation whose explicit discretization may be done directly by
the following way.

Denoting for  $x\in m_i \,,\, f_{\alpha}(t_n,x,\xi)~=~\mathcal{M}_{\alpha,i}^n(\xi) = \mathcal{M}_{\alpha}(A_{\alpha,i}^n,Q_{\alpha,i}^n,\xi)$
the Maxwellian state associated to  $A_{\alpha,i}^n$ and $Q_{\alpha,i}^n$, a finite volume discretization of Equation \eqref{eqcin} leads to:
\begin{equation} \label{cindiscret}
f_{\alpha,i}^{n+1}(\xi) = \mathcal{M}_{\alpha,i}^n(\xi)+\frac{\Delta t}{h_i}\,
\xi \, \left(\mathcal{M}_{\alpha,i+\frac{1}{2}}^-(\xi)-\mathcal{M}_{\alpha,i-\frac{1}{2}}^+(\xi)\right)
\end{equation}
where the fluxes $\mathcal{M}_{\alpha,i+\frac{1}{2}}^\pm$ have to take into account the discontinuity of
the source term $\phi_{\alpha}$ at the cell interface $x_{i+1/2}$. This is the principle of interfacial  source upwind. 
Indeed, noticing that the fluxes can also be written as:
\begin{equation*}
\mathcal{M}_{\alpha,i+\frac{1}{2}}^-(\xi) = \mathcal{M}_{\alpha,i+\frac{1}{2}} +
\left(\mathcal{M}_{\alpha,i+\frac{1}{2}}^- - \mathcal{M}_{\alpha,i+\frac{1}{2}}\right)
\end{equation*}
the quantity $\delta \mathcal{M}_{\alpha,i+\frac{1}{2}}^- = \mathcal{M}_{\alpha,i+\frac{1}{2}}^- -
\mathcal{M}_{\alpha,i+\frac{1}{2}}$
holds for the discrete contribution of the source term $\phi_{\alpha}$ in the system for negative
velocities $\xi \leq 0$  due to the upwinding of the source term.
Thus $\delta \mathcal{M}_{\alpha,i+\frac{1}{2}}^-$ has to vanish for positive velocity $\xi > 0$,
as proposed by the choice of the interface fluxes below.
Let us now detail our choice for the fluxes $\mathcal{M}_{\alpha,i+\frac{1}{2}}^\pm$ at the interface. It
can be justified by using a generalized characteristic method for Equation \eqref{KineticFormulationPFS}
(without the collision kernel) but we give instead a presentation based on some physical energetic
balance. The details of the construction of these fluxes by the general characteristics  method (see \cite[Definition 2.1]{Da89}) is
done in \cite[Chapter 2]{TheseErsoy}.

In order to take into account the neighboring cells by means of a natural interpretation of the
microscopic features of the system, we formulate a peculiar  discretization for the fluxes in
\eqref{cindiscret}, computed by the following upwinded formulas:
\begin{equation}\label{MicroscopicInterfaceFluxes}
\begin{array}{lll}
\mathcal M_{\alpha,i+1/2}^{-}(\xi) &=&  
\overbrace{\dsp \mathds{1}_{\{\xi>0\}}\mathcal M_{\alpha,i}^n(\xi)}^{{\textrm{positive  transmission}}}
+\overbrace{\mathds{1}_{\{\xi<0,\xi^2-2g\Delta\phi^n_{\alpha,i+1/2}<0\}}\mathcal M_{\alpha,i}^n(-\xi)}^{\textrm{reflection}}\\ 
&+& \underbrace{\dsp \mathds{1}_{\{\xi<0,\xi^2-2g\Delta\phi^n_{\alpha,i+1/2}>0\}}
\mathcal M_{\alpha,i+1}^n\left(\dsp-\sqrt{\xi^2-2g\Delta\phi^n_{\alpha,i+1/2}}\right)}_{\textrm{negative   transmission}},\\
 & & \\
\mathcal M_{\alpha,i+1/2}^{+}(\xi) &=& \overbrace{\dsp \mathds{1}_{\{\xi<0\}}\mathcal M_{\alpha,i+1}^n(\xi)}^{\textrm{negative transmission}}
+\overbrace{\mathds{1}_{\{\xi>0,\xi^2+2g\Delta\phi^n_{\alpha,i+1/2}<0\}}\mathcal M_{\alpha,i+1}^n(-\xi)}^{\textrm{reflection}}\\ 
&+& \underbrace{ \dsp \mathds{1}_{\{\xi>0,\xi^2+2g\Delta\phi^n_{\alpha,i+1/2}>0\}}
\mathcal M_{\alpha,i}^n\left(\dsp\sqrt{\xi^2+2g\Delta\phi^n_{\alpha,i+1/2}}\right)}_{\textrm{positive  transmission}}.
\end{array}\,
\end{equation}
The term $\Delta\phi^n_{\alpha,i\pm 1/2}$ in \eqref{MicroscopicInterfaceFluxes}  is the upwinded source term
\eqref{PFSSourceTermPhi}. It also plays the role of the potential barrier:  
the term $\xi^2\pm 2g\Delta\phi^n_{\alpha,i+1/2}$ is the jump condition for a particle with a kinetic speed $\xi$  
which is necessary  to
\begin{itemize}
\item be reflected: this means that the particle has not enough kinetic energy $\xi^2/2$ 
to overpass the potential barrier (reflection term  in \eqref{MicroscopicInterfaceFluxes}),
\item overpass the potential barrier with a positive speed 
(positive  transmission term in \eqref{MicroscopicInterfaceFluxes}),
\item overpass the potential barrier with a negative speed 
(negative  transmission term in \eqref{MicroscopicInterfaceFluxes}).
\end{itemize}
Taking an approximation of the non-conservative product $\phi_{\alpha}$  defined by Equation
\eqref{DefinitionOfNonConservativeProductPhi}, 
the potential barrier $\Delta \phi^n_{\alpha,i+ 1/2}$ has the following expression:
\begin{equation*}
\Delta\phi_{\alpha,i+1/2}^n= \jump{\WW_{\alpha}}(t_n)\cdot  \BB_{\alpha}\left(t_n,\Psi\left(\frac{1}{2},\WW_{\alpha,i}(t_n),
\WW_{\alpha,i+1}(t_n)\right)\right)
\end{equation*}
which is an approximation of $\Delta\phi_{\alpha}$ by the midpoint quadrature formula. 

Since we neglected the collision term, it is clear that $f_{\alpha}^{n+1}$ computed by the discretized
kinetic equation \eqref{cindiscret} is no more a Gibbs equilibrium. Therefore, to recover the macroscopic variables $A_{\alpha}$ and $Q_{\alpha}$,
according to  the identities \eqref{macroA}-\eqref{macroQ}, we set:
\begin{equation*}
 \UU_{\alpha,i}^{n+1}=
\left(\begin{array}{l}
A_{\alpha,i}^{n+1} \\
Q_{\alpha,i}^{n+1}
\end{array}
\right)
\overset{def}{=} \int_\R \left(\begin{array}{l}
1 \\
\xi
\end{array}
\right) f_{\alpha,i}^{n+1} \,d\xi
\end{equation*}
In fact at each time step, we projected $f_{\alpha}^{n}(\xi)$ on $\mathcal{M}_{\alpha,i}^n(\xi)$, 
which is a way to perform all collisions at once and to recover a Gibbs equilibrium
without computing it.

Now, we can integrate the discretized kinetic equation \eqref{cindiscret} against 1 and $\xi$ to obtain the macroscopic kinetic scheme:
\begin{equation} \label{macrodiscret}
\UU_{\alpha,i}^{n+1} = \UU_{\alpha,i}^n +\frac{\Delta t^{n}}{h_i}\left(\F_{\alpha,i+\frac{1}{2}}^- - \F_{\alpha,i-\frac{1}{2}}^+ \right) \: .
\end{equation}
The numerical fluxes are thus defined by the kinetic fluxes as follows:
\begin{equation}\label{fluxdiscret}
\F_{\alpha,i+\frac{1}{2}}^\pm
\overset{def}{=} \int_\R \xi \left(\begin{array}{l}
1 \\
\xi
\end{array}
\right) \mathcal{M}^\pm_{\alpha,i+\frac{1}{2}}(\xi)\,d\xi
\end{equation}
Computing the macroscopic state $\UU_{\alpha}$ by Equations \eqref{macrodiscret}-\eqref{fluxdiscret} is not easy if the function $\chi$ 
verifying the properties \eqref{propchi} is not compactly supported (see \cite{PS01}). 

Moreover, as we shall see in the next proposition, a CFL conditions is needed to obtain a scheme that preserves the positivity of the wet area.
\begin{proposition}\label{SchemaCinProprieteClassique}
Let $\chi$ be a compactly supported function verifying  \eqref{propchi} and note $[-K,K]$ its support. 
The kinetic scheme \eqref{macrodiscret}-\eqref{fluxdiscret} has the following properties:
\begin{enumerate}
\item The kinetic scheme is a $A_{\alpha}$ conservative  scheme,
\item Assume the following CFL condition  
\begin{equation}\label{CFLcondition}
\Delta  t^n \max_{\alpha,i}\left(\abs{u_{\alpha,i}^n}+ K\,b_{\alpha,i}^n\right)\leqslant \min_i h_i
\end{equation} holds. 
Then the kinetic scheme keeps the wet area $A_{\alpha}$ positive i.e:
$$\mbox{ if, for every } i\in [0,N+1] \,,\, A^{0}_{\alpha,i} \geqslant 0 \mbox{ then, for every } i\in [0,N+1] \,,\, A_{\alpha,i}^n \geqslant 0.$$
\item The kinetic scheme treats ``naturally'' flooding and drying zones for the water layer, and the apparition of vacuum for the air layer.
\end{enumerate}
\end{proposition}
\begin{proof}
We will adapt the proof of \cite{PS01} to show the three properties that verify  the kinetic scheme.

1. Let us denote the first component of the discrete fluxes \eqref{fluxdiscret}
$\left(F_{A_{\alpha}}\right)^\pm_{i+\frac{1}{2}}$:
\begin{equation*}
\left(F_{A_{\alpha}}\right)^\pm_{i+\frac{1}{2}}
\overset{def}{=} \int_\R \xi \mathcal{M}^\pm_{\alpha,i+\frac{1}{2}}(\xi)\,d\xi
\end{equation*}
An easy computation, using the change of variable $\mu = \abs{\xi}^2 - 2g\Delta\phi_{\alpha,i+\frac{1}{2}}^{n}\,$, allows us
to show that:
\begin{equation*}
\left(F_{A_{\alpha}}\right)^+_{i+\frac{1}{2}} = \left(F_{A_{\alpha}}\right)^-_{i+\frac{1}{2}} \: .
\end{equation*}
2. Suppose that for every $i\in [0,N+1]$, $A_{\alpha,i}^{n}>0$. Let us note $\xi_{\pm} = \max(0,\pm\xi)$  and  $\sigma=\dsp\frac{\Delta t^n}{\min_{i} h_i}$.
From Equation \eqref{cindiscret},  we get the following inequalities:
$$\begin{array}{lll}
 f_{\alpha,i}^{n+1}(\xi)  & \geqslant&  (1-\sigma |\xi|)\mathcal M_{\alpha,i}^n(\xi)\\
                 &   & + \sigma \xi_+ \bigg[\mathds{1}_{\{\xi^2+2g\Delta \phi_{\alpha,i+1/2}<0\}}\mathcal
M_{\alpha,i}^n(-\xi)\\ 
                 &   & +  \mathds{1}_{\{\xi^2+2g\Delta \phi_{\alpha,i-1/2}>0\}}\mathcal
M_{\alpha,i-1}^n\left(\displaystyle\sqrt{\xi^2+2g\Delta
\phi_{\alpha,i+1/2}}\right)\bigg]  \\
		 &   &  +  \sigma \xi_- \bigg(\mathds{1}_{\{\xi^2-2g\Delta \phi_{\alpha,i+1/2}<0\}}\mathcal
M_{\alpha,i}^n(-\xi)\\  
                 &   &  + \mathds{1}_{\{\xi^2-2g\Delta \phi_{\alpha,i-1/2}>0\}}\mathcal
M_{\alpha,i+1}^n\left(\displaystyle-\sqrt{\xi^2-2g\Delta
\phi_{\alpha,i+1/2}} 
\right)\bigg)  \ . 
\end{array}$$
Since the function $\chi$ is compactly supported $\textrm{ if } \abs{\xi - u_{\alpha,i}^n}\geqslant K b_{\alpha,i}^n \textrm{ then } \mathcal{M}_{\alpha,i}^n(\xi) = 0 .$
Thus $$f_{\alpha,i}^{n+1}(\xi)\geqslant 0 \textrm{ if } \abs{\xi - u_{\alpha,i}^n}\geqslant  K b_{\alpha,i}^n \,,$$ 
as a sum of non negative terms.\\
On the other hand, for $\dsp \abs{\xi - u_{\alpha,i}^n}\leqslant K b_{\alpha,i}^n$, using the CFL condition   $0 < \sigma \abs{\xi}\leqslant 1$, 
for all $i$, $f_{\alpha,i}^{n+1}\geqslant 0$ since it is a convex combination of non negative terms. 

Finally we have $\forall i \in [0,N+1] \,,\, f^{n}_{\alpha,i}  \geq 0$. Since $ \forall i \in [0,N+1] \,,\, A_{\alpha,i}^{n+1} = 
\displaystyle \int_{\R} f_{\alpha,i}^{n+1}(\xi)\,d\xi,$ we finally get  
$ \forall i \in [0,N+1] \,,\ A_{\alpha,i}^{n+1} \geq 0.$\\
$3.$ Suppose $A_{i}^{n} = 0$. Using the definition of $\mathcal{M}_w$, and the fact that the function $\chi$ is compactly
supported, the only term that may cause problem is $\dsp \frac{A}{b_w(t,x)}.$
But since 
$$\dsp \frac{A}{b(t,x)} \sim \sqrt{\frac{A}{g I_1(x,A)\cos\theta}} \mbox{ when }A \sim 0,$$ 
we get 
$$\dsp  \lim_{\underset {A \geq 0} {A \rightarrow 0}} \frac{A}{b_w(t,x)}  =  0 \; ,$$ 
thus $\mathcal{M}_{w,i}^n(\xi) = 0$. This is the reason why we say
that the kinetic scheme treats ``naturally'' the drying and flooding zones.
In the same way, for the air layer, the only term which may pose problem is $\frac{M}{b_a}$ which is, when $M$ goes to $0$, equivalent to $M^{5/2-\gamma}$ 
where $5/2-\gamma=11/10>1$. Therefore, we say that the numerical scheme computes ``naturally'' the apparition/vanishing of vacuum.
\end{proof}
Let us finally point out that the $\chi$ function could be chosen undifferently.

\subsection{The ``decoupled boundary conditions''}\label{SectionBoundaryTreatment}
Without loss of generality, we consider in what follows the case of a uniform section with $\theta = 0$ slope.
Let us recall that $x_{1/2}$ and $x_{N+1/2}$ are respectively the upstream and the downstream ends of the pipe. 
At this stage, we have computed all the ``interior'' states at time $t^{n+1}$, that is $(\WW_{i}^{n+1})_{i=1,N}:=
(M_{i}^{n+1},D_{i}^{n+1},A_{i}^{n+1},Q_{i}^{n+1})_{i=1,N}$ are computed.

The upstream state $\WW_{0}$ corresponds to the mean value of $(M,D,A,Q)$ on the ``fictitious'' 
cell $m_0=(x_{-1/2},x_{1/2})$ (at the left of the upstream boundary of the pipe) 
and the downstream state $\WW_{N+1}$ corresponds to the mean value of $(M,D,A,Q)$ on the ``fictitious'' cell $m_{N+1}=(x_{N+1/2},x_{N+3/2})$ 
(at the right of the downstream end of the pipe).

Usually,  we have to prescribe at least two boundary conditions related to the state vectors
$\WW_0^n$ and $\WW_{N+1}^n$. There is generically, that is to say for supercritical flows,  ``two incoming characteristics curve'' for the upstream and  two for 
outgoing characteristics for the downstream when the eigenstructure is known which is clearly not the case in the full System \eqref{BilayerModel}. 

In order to find a complete state for the upstream boundary, $\WW_{0}^{n+1}$, and the downstream boundary, $\WW_{N+1}^{n+1}$, in a consistent way, 
following the recent work of Bourdarias \emph{et al.} \cite{BEG11_2}, one can use the so-called ``kinetic boundary conditions'' 
for which the eigenvalues are not required 
and can be applied straightforwardly (at the present stage, such boundary conditions are not yet implemented).
We present instead a new approach, based on the solution of a linearized Riemann problem with source term, which introduces the concept of 
``decoupled boundary conditions''. Indeed, the full eigenstructure is not available for the full system but dealing with each layers separatively, 
we have two systems which write:
\begin{equation}\label{eau}
\partial_t \UU + \mathcal{D}_{w}(\WW) \partial_x \UU = \left(\begin{array}{c} 0\\a c^2_a \partial_x M  \end{array}\right)
\end{equation}
and
\begin{equation}\label{air}
\partial_t \VV + \mathcal{D}_{a}(\WW) \partial_x \VV = \left(\begin{array}{c} 0\\- c^2_a \frac{A M}{(S-A)^2} \partial_x A \end{array}\right)
\end{equation}
where $\UU = (A,Q)^t$, $\VV = (M,D)^t$, 
$\mathcal{D}_{w}(\WW) = \left(\begin{array}{cc} 0 & 1 \\ c_m^2-u^2 & 2 u \end{array}\right)$ and 
$\mathcal{D}_{a}(\WW) = \left(\begin{array}{cc} 0 & 1 \\ c_a^2-v^2 & 2 v \end{array}\right)$ 
for which eigenstructure is well known and reads as follows for the water layer 
$$\lambda_w^1(\WW) = u-c_m,\quad \lambda_w^2(\WW) = u+c_m\, ,$$ 
and for the air layer: 
$$\lambda_a^1(\WW) = v-c_a,\quad \lambda_a^2(\WW) = v+c_a\, .$$ 
In theses equations, $u= Q/A$ and $v=D/M$ stands for the mean water and air speed.

Thus, one can apply the usual procedure \cite[Chapter 2]{TheseErsoy} (see, also \cite{BEG11_2}) which consists to assume that the inflow is subcritical, i.e. 
we have to prescribe at least two boundary conditions one for each layers corresponding to the 
``two incoming decoupled characteristics curve'' for the upstream and  two for 
outgoing characteristics for the downstream. 
Therefore, for instance, at the upstream end of the pipe, one of the following boundary conditions may be prescribed  for the water layer (as usual) 
(we omit the index $n+1$ for the sake of clarity): 
\begin{enumerate}
 \item the water level is prescribed. So let $H_{up}(t)$ be a given function of time. Then we have:  
\begin{equation}\label{amonthauteur}
\forall t > 0 ,\quad  \dsp \mathcal{H}(A_0(t))\cos\theta_0+\mathbf{Z_0} = H_{up}(t) \ .
\end{equation}
 \item the discharge is prescribed. So let $Q_{up}(t)$ be a given function of time. Then we have: 
 \begin{equation}\label{UpstreamDischarge}
\forall t > 0 ,\quad \dsp Q_0(t) = Q_{up}(t) \ .
\end{equation}
\end{enumerate}
Similarly for the air layer, to be consistent with respect to  the previous boundary conditions, we prescribe the air density at the upstream and downstream end.
Thus, for instance, 
if $(\rho_{up}(t),A_{up}(t))$ are given, 
by definition $M_{up}(t):=\frac{\rho_{up}(t)}{\rho_0} (S-A_{up}(t))$ is given and we have to compute $D_{up}$ and $Q_{up}$. 
In the same way, if $(\rho_{up}(t),Q_{up}(t))$ are given then we have to compute $D_{up}$ and $A_{up}$. Once $A_{up}(t)$ is computed, one has 
$M_{up}(t):=\frac{\rho_{up}(t)}{\rho_0} (S-A_{up}(t))$  by definition. 
For each cases, we have to define two equations, which can be obtained by the usual procedure which consists to solve two  
linearized Riemann problems for each layers, which can be written into a generic form as follows: 
\begin{equation}\label{mixte}
\left\{
\begin{array}{l}
\partial_t \UU_{\alpha} + \widetilde{\mathcal{D}_{w}} \partial_x \UU_{\alpha} = G_{\alpha}(\WW) \\
\UU_{\alpha}(x,t_n) =  
\left\{
\begin{array}{lll}
\UU_{\alpha,0}^n \textrm{ if } x<0,\\
\UU_{\alpha,1}^n \textrm{ if } x>0
\end{array}
\right.
\end{array}
\right.
\end{equation}
where $\widetilde{\mathcal{D}_{\alpha}}$ stands for the constant matrix $\mathcal{D}_{\alpha}(\widetilde \WW)$ with 
$\widetilde\WW = \frac{\WW_0^n+\WW_1^n}{2}$, $G_{\alpha}(\WW)$ the source term appearing in Equation \eqref{eau} or \eqref{air} where 
 $\alpha$  stands for the subscript $a$ or $w$. 
The previous Riemann problem allows to define the missing two equations required to close the system of two unknowns according to the prescribed conditions, and 
we deduce the two missing equations which read:
\begin{equation}\label{BCeau}
(Q_0^{n+1}-Q_1^{n})-(A_0^{n+1}-A_1^{n})(\widetilde{u}+\widetilde{c_m}) = 
\frac{\Delta t}{2} \left(\left(c_a^2 \frac{A}{S-1}\right)_0^{n+1}+\left(c_a^2 \frac{A}{S-1}\right)_1^{n}\right) (M_0^{n+1}-M_1^{n})
\end{equation}
and 
\begin{equation}\label{BCair}
(D_0^{n+1}-D_1^{n})-(M_0^{n+1}-M_1^{n})(\widetilde{v}+\widetilde{c_a}) = 
\frac{\Delta t}{2} \left(\left(c^2_a \frac{M}{S-A}\right)_0^{n+1}+\left(c^2_a \frac{M}{S-A}\right)_1^{n}\right) (A_0^{n+1}-A_1^{n})
\end{equation}
Let us remark that, when $(\rho_{up}(t),A_{up}(t))$ are prescribed  $Q_0^{n+1}$ and $D_0^{n+1}$ are respectively explicitly provided by Equation 
\eqref{BCeau} and \eqref{BCair}. By the contrary, if $(\rho_{up}(t),Q_{up}(t))$ are prescribed, in order to get  
$A_0^{n+1}$ and $D_0^{n+1}$, we solve the non linear system formed by Equation \eqref{BCeau} and \eqref{BCair}. 

Since the subcritical flows is the only one for which we can write ``real boundary conditions'', we process by \emph{a priori}, i.e. after the computation 
of the subcritical flows as described by the above-mentioned process,  we check the validity by computing the sound speed.  
If such assumption is incorrect, we then switch to the two following configurations:
\begin{itemize}
 \item an incoming supcritical flows is observed then the critical flows is prescribed (to complete the missing equation),
 \item an outgoing supcritical flows is observed then we prescribe a free boundary.
\end{itemize}
We apply such \emph{a priori} process to deal with such inflow and outflow conditions. 

Finally, let us point out that, in previous works the authors have shown its efficiency by validating the treatment of boundary conditions through 
different test cases in the context of transient mixed water flows in closed water pipes (see for instance \cite{BEG11_1,BEG11_2}).

\section{Numerical tests}\label{Numeric}
The numerical experiments  consist in studying the air entrainment effect in a closed pipe. 
Since experimental data for such  flows in any pipes are not available, we focus only on  the qualitative behavior of the proposed
model and the presented method for different upstream and downstream conditions. 

For all tests, we assume that the air  is injected inside the pipe under consideration in order to ensure a constant density of air at upstream 
and downstream end. For this purpose, we set the upstream and downstream air density $\rho$ to $\rho_a$: 
we have chosen $\rho_{a} = 1.29349\, 10^{-2}\:  kg/m^3$ and $\rho_{a} = 1.29349\:  kg/m^3$ and  as an initial air density, we always set 
$\rho_{init}(x) = \rho_a$.  
Then, the 
unknown state vector $(A,Q,M,D)$ at the upstream and downstream end evolves as described in Section \ref{SectionBoundaryTreatment}. 
We also  represent the non hyperbolic region.  
Finally, at the end of the section, we present a qualitative experiment with respect to the sensitivity of the 
numerical study with respect to the spatial mesh size  when the system is mainly  non hyperbolic.

For all tests, the function $\chi$ satisfying \eqref{propchi}
is defined by:
\begin{equation*}\label{chifunction}
\chi=\frac{1}{2\sqrt{3}}\mathds{1}_{[-\sqrt{3},\sqrt{3}]} \ .
\end{equation*}
\subsection{Non constant upstream boundary conditions}
We consider an horizontal  pipe of circular cross-section of $2 \: m^2$ is $10\: m $ long. 
The altitude of the upstream end of the pipe is $0 \: m$.  
We start the numerical simulation with a  steady state with the water height $$h_{init}(x) = 0.05\:  m$$ 
(i.e. $A_{init}(x) = 2.0\, 10^{-2}\: m^2$), $M_{init}(x):=\frac{\rho_{a}}{\rho_0}(S(x)-A_{init}(x)) = 4.0\,  10^{-3}\: m^2$, 
the water discharge $Q_{init}(x) = 0\: m^3/s$, the air discharge $$D_{init}(x) = 0\: m^3/s\, ,$$ for all $x\in [0,10]$.   
The upstream water height is set to $$h_0(t) = 0.05(1+ t) \mathds{1}_{[0,5]}(t) + 0.3  \mathds{1}_{(5,50]}(t)$$ while 
the downstream water discharge is kept constant equal to $$D_{N+1}(t)=  0\: m^3/s\, .$$ 

\begin{figure}[H]
\subfigure[Prescribed water height \label{CompUpstreamBC}]
{
\includegraphics[height = 5.8cm]{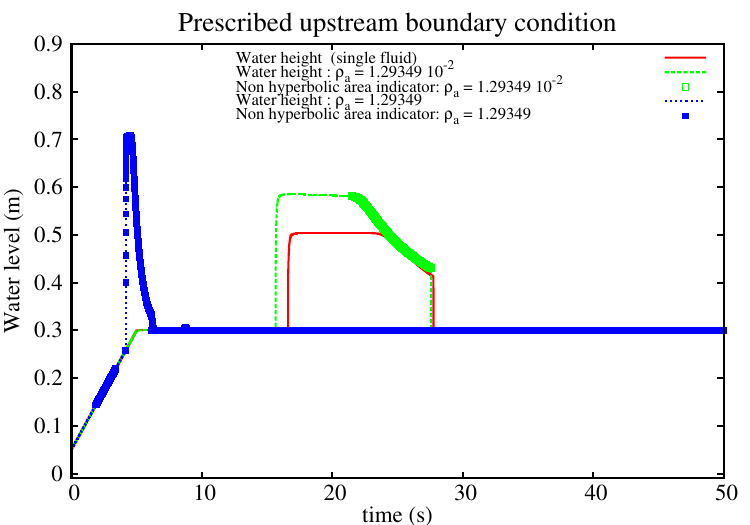}
} 

\subfigure[Water speed \label{CompUpstreamWaterSpeed}]
{
\includegraphics[height = 5.8cm]{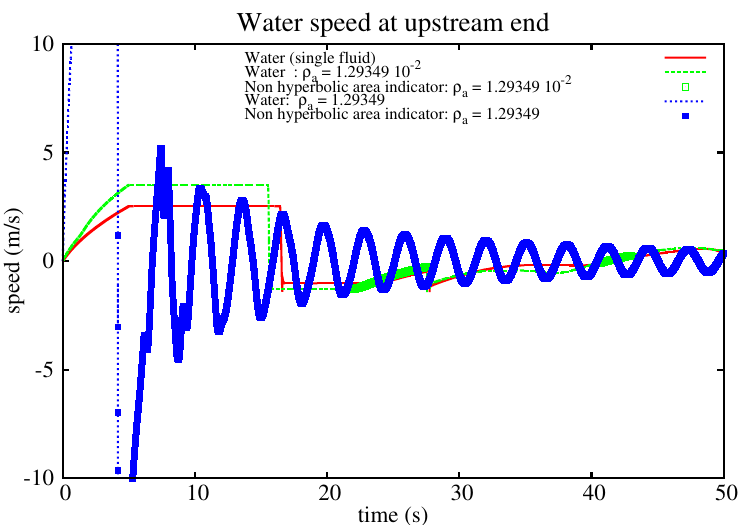}
} 

\subfigure[Air speed \label{CompUpstreamAirSpeed}]
{
\includegraphics[height = 5.8cm]{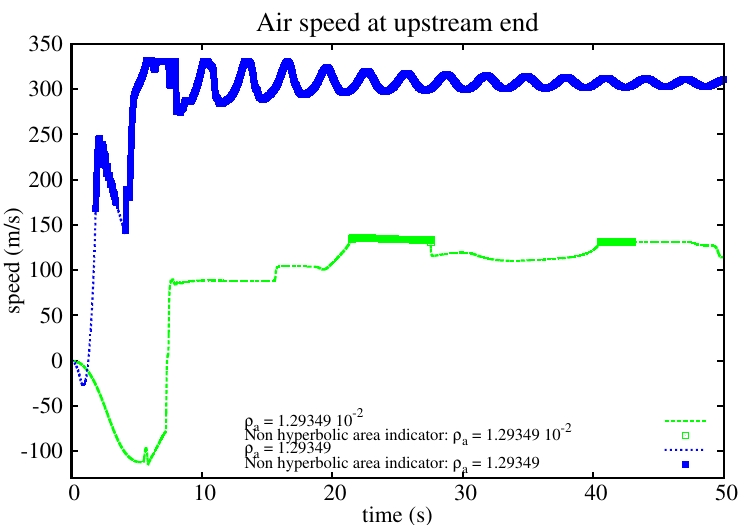}
}
\caption{Non constant upstream boundary conditions: numerical results at upstream end\label{CompUpstream}}
\end{figure}

\begin{figure}[H]
\subfigure[Water height \label{Comp5Piezo}]
{
\includegraphics[height = 5.8cm]{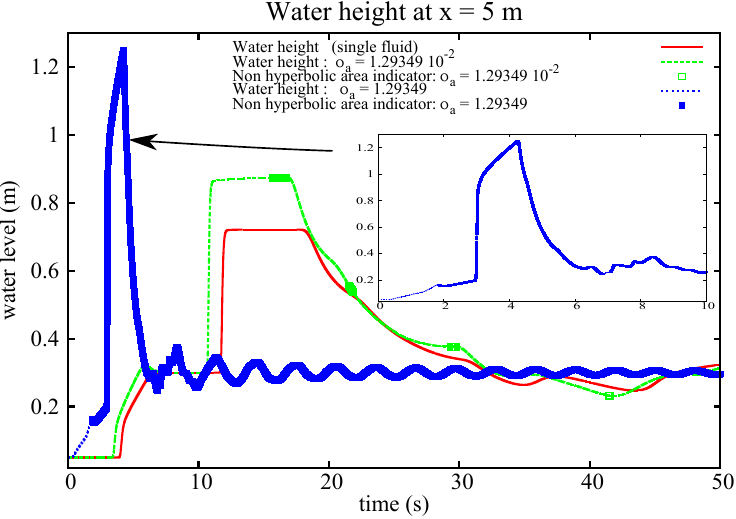}
} 

\subfigure[Water speed \label{Comp5WaterSpeed}]
{
\includegraphics[height = 5.8cm]{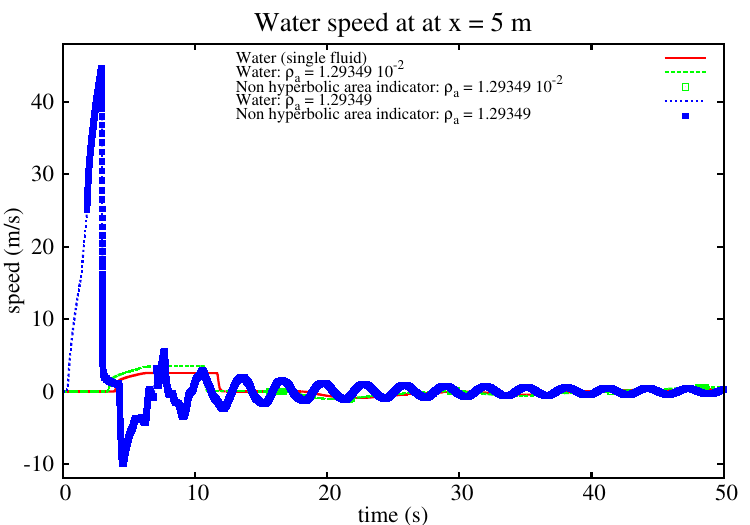}
} 

\subfigure[Air speed \label{Comp5AirSpeed}]
{
\includegraphics[height = 5.8cm]{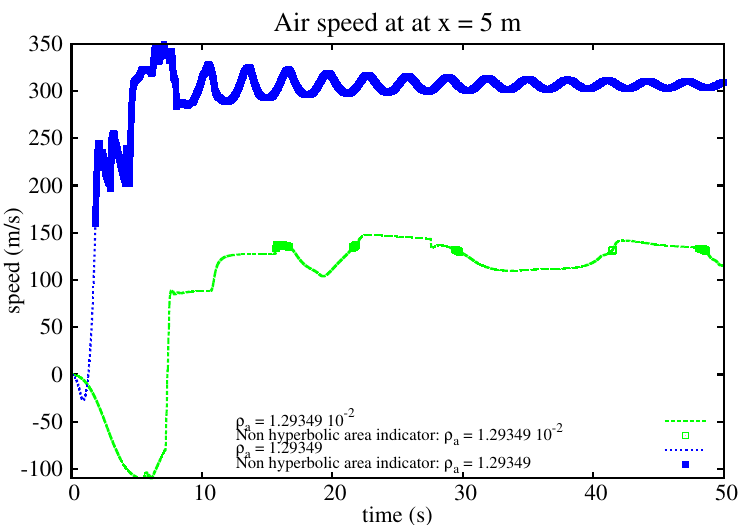}
}
\caption{Non constant upstream boundary conditions: numerical results at $x=5\: m$\label{Comp5}}
\end{figure}

\begin{figure}[H]
\subfigure[Water height \label{CompDownstreamPiezo}]
{
\includegraphics[height = 5.8cm]{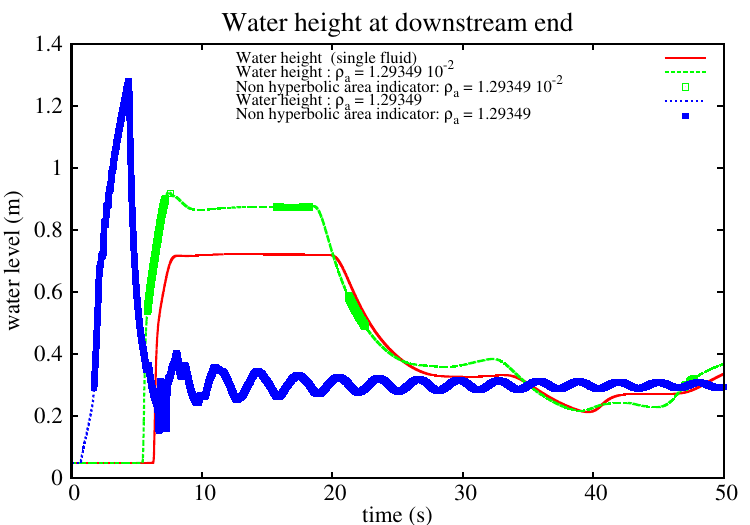}
}

\subfigure[Prescribed downstream discharge\label{CompDownstreamBC}]
{
\includegraphics[height = 5.8cm]{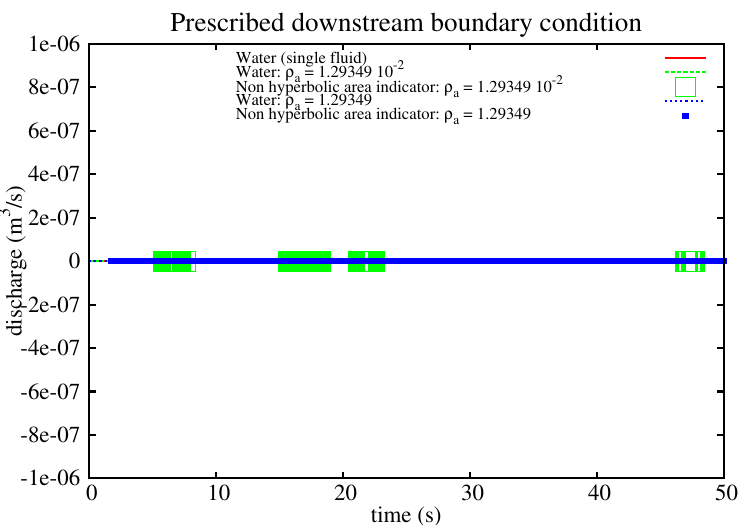}
} 

\subfigure[Air speed \label{CompDownstreamAirSpeed}]
{
\includegraphics[height = 5.8cm]{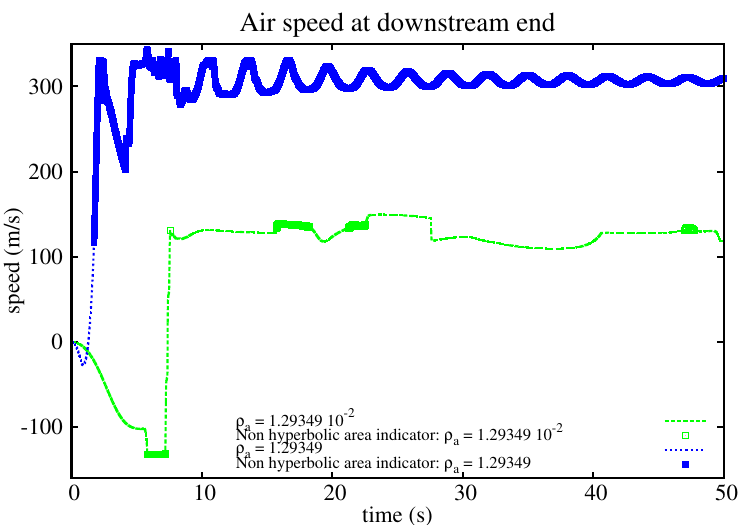}
}
\caption{Non constant upstream boundary conditions: numerical results at downstream end\label{CompDownstream}}
\end{figure}

These boundary conditions are chosen such that the air-layer model will lose its hyperbolicity for the three following numerical experiments:
\begin{enumerate}
 \item the free surface system (called thereafter ``single fluid'').
 \item the two-layer system with $\rho_{a} = 1.29349\, 10^{-2}\:  kg/m^3$ fixed at upstream and downstream end.
 \item the two-layer system with $\rho_{a} = 1.29349\:  kg/m^3$ fixed at upstream and downstream end.
\end{enumerate}
The single fluid case is solved by a  kinetic scheme   \cite{PS01,BEG11_1,BEG11_2} while the others tests by the presented 
``two-layer kinetic scheme''. The obtained results are  compared 
on \resim\ref{CompUpstream}, \ref{Comp5} and \ref{CompDownstream}  at upstream end, at $x=5\: m$ and at downstream end. 
We also represents the non hyperbolic region  by square points.  

Other parameters are: 
$$
\begin{array}{lcl}
\hbox{Discretisation points}& : & 100,\\ 
\hbox{Delta x }(m)& : &0.1,\\
\hbox{CFL }& : & 0.95,\\
\hbox{Simulation time } (s) & : &50.
\end{array}
$$

The case of the so-called ``single fluid'' has been recently validated through 
several numerical experiments (see  \cite{BEG11_1,BEG11_2}). 
Since we use the same methodology, the two-layer kinetic scheme presented here will be compared to the ``single fluid'' case to show
the effect of air entrainment. 
The presented first numerical experiments simulate the propagation of an hydraulic bore emerging from the 
upstream end as we  increase the upstream water height 
(see \resim\ref{CompUpstreamBC} for $t\leqslant 5\: s$).  Then, it reaches the downstream end at, approximatively, 
$t = 7\: s$ the first time (see \resim\ref{CompDownstreamPiezo}) and it is reflected since the dowstream discharge is kept constant equal to zero. 
Thereafter, it comes back to the upstream end at, approximatively, $t = 17\: s$ and until $t\approx 27\: s$ a free boundary conditions is 
prescribed (which corresponds to an outgoing subcritical flows, i.e. $u_0<-c_0$, see \resim\ref{CompUpstreamBC} and 
Section \ref{SectionBoundaryTreatment} for further details). This wave propagates alternatively from 
upstream to downstream before to be damped which is observed on \resim\ref{Comp5Piezo} and \ref{CompDownstreamPiezo} 
for $t\geqslant 17\: s$.\newline

Now, let us focus on the behavior of the numerical solution when the air is taken into account and let's compare the numerical solution with the one obtained for the ``single fluid''. 

As a first remark, let us  point out that the shape of the numerical solution of the two-layer system is similar up to 
a non linear transformation: the shape is rescaled with a higher amplitude. It means that the 
propagation of the hydraulic bore is also observed when air is taken into account but with a greater celerity: this means that the hydrodynamic 
behavior of the water is slightly and sharply modified (following the value of $\rho_a$) under the interaction of air. 
Globally, on \resim\ref{CompUpstream}, \ref{Comp5} and \ref{CompDownstream}, this fact is clearly visible 
when $\rho_{a} = 1.29349\, 10^{-2}\:  kg/m^3$ is prescribed  
while we have to pay our attention to the first seconds when $\rho_{a} = 1.29349\:  kg/m^3$. For instance, on 
\resim\ref{Comp5WaterSpeed}, the shape for $t\leqslant 30\: s$ of the numerical solution of the single fluid has 
to be compared with the shape for $t\leqslant 8\: s$ in this case.

As a second remark, let us point that the celerity of this phenomenon is an increasing function of $\rho_a$. For instance,  
when $\rho_{a} = 1.29349\:  kg/m^3$, the hydraulic bore reaches immediately the downstream end with a higher amplitude 
than the case when $\rho_{a} = 1.29349\, 10^{-2}\:  kg/m^3$ is prescribed which approximatively requires $7\: s$ 
(see for example \resim\ref{CompDownstreamPiezo}). These phenomena are the consequence of the prescribed boundary conditions and above all 
of the compressibility of the air layer which  interacts with the water layer both confined in a closed area. As the effect are stronger in the case 
(2) than in the case (3), from now on we focus only in the test case $(3)$.

At time $t=0\: s$, the air layer as well as the  water layer are both in a closed environment. 
Between, $0\leqslant t \leqslant 1$, increasing the water height increases the water speed and the air  escapes from the upstream (i.e. with a negative speed) 
as well as in the whole pipe (see \resim\ref{CompUpstreamAirSpeed}, \ref{Comp5AirSpeed} and \ref{CompDownstreamAirSpeed}) 
due to  the prescribed boundary conditions. Indeed, we increase the water level at upstream end and we set the water discharge to be $0$ at downstream end. 
Therefore, the air necessary escapes throughout the upstream end.  
After this time period, see \resim\ref{CompUpstreamWaterSpeed} and \ref{CompUpstreamAirSpeed} for  $1\leqslant t \leqslant 5$, 
the air speed reaches a minimum value at $t\approx1\: s$ and since 
we are forcing the air density to be a constant, the upstream boundary equation \eqref{BCair} ensures that the upstream air incomes is enough to keep it constant.
Consequently, the water speed is again accelerated and  we observe a maximum value at $t\approx 5\: s$. 
An the same time, as the water flow is accelerated, the hydraulic bore  at $t\approx 5\: s$, after being reflected at the downstream end, is 
present at the upstream end with a negative speed. This explains why the water speed becomes negative.  
These phenomena is then reproduce with attenuation of the amplitude of the wave as represented by the oscillation on 
\resim\ref{CompUpstreamWaterSpeed} and \ref{CompUpstreamAirSpeed}. The same explanation holds at $x=\: 5m$ 
and at the downstream end as we observe on \resim\ref{Comp5WaterSpeed}, \ref{Comp5AirSpeed} and  \resim\ref{CompUpstreamWaterSpeed}, \ref{CompUpstreamAirSpeed}.

These phenomena are reproduced many times, as far as, the initial waves propagates from the upstream to the downstream end until it vanishes due to 
some dissipative effects. These corresponds to the oscillations observed on \resim\ref{CompUpstream}, \ref{Comp5} and \ref{CompDownstream} for the test 
case $(3)$.  \newline

From now on, let us focus on the non hyperbolic aspect and let us 
emphasize that the two-layer system loses its hyperbolicity in both tests $(2)$ and $(3)$ and  
the area of the non hyperbolic region is  increasing with respect to $\rho_a$ as we can observe on \resim\ref{CompUpstream}, \ref{Comp5}, 
\ref{CompDownstream} (which numerically confirms  the analysis done in Section 
\ref{SubsectionBilayerModel}, see Theorem \ref{BilayerModelThm}).

\emph{Let us point out that the oscillatory phenomenon is not the so-called geometric growing instabilities} and thus not at all related  
to the loss of hyperbolicity of the two-layer system.  Indeed, on one hand, such oscillations are also present for the ``single fluid'' case 
but with very small amplitude as we can observe on \resim\ref{CompDownstreamPiezo}. 
On the other hand, in view of the previous explanations of the origin of the oscillations, these phenomena have to be compared 
with the so-called \emph{waterhammer} occurring during the filling of the pipe (see,  for instance \cite{BEG09_1,BEG11_2}). 

These phenomena seems to be a ``airhammer'' for which the number of oscillations increases with the density of the air (as in the case 
of pressurized flows with  the sound speed, see \cite{BEG09_1,BEG11_2}).  This is thus a dynamic consequence of the compressibility of the air layer 
%
\subsection{Non constant downstream water discharge}
\begin{figure}[H]
\subfigure[Water height \label{Bar5Piezo}]
{
\includegraphics[height = 5.8cm]{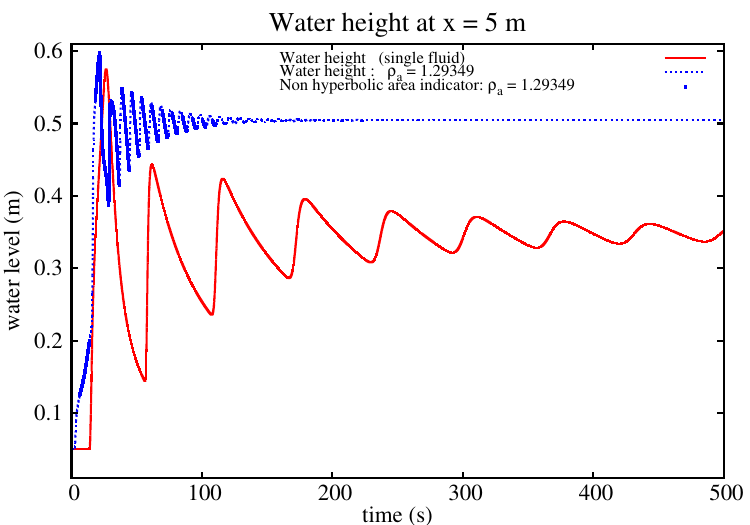}
} 

\subfigure[Water speed \label{Bar5WaterSpeed}]
{
\includegraphics[height = 5.8cm]{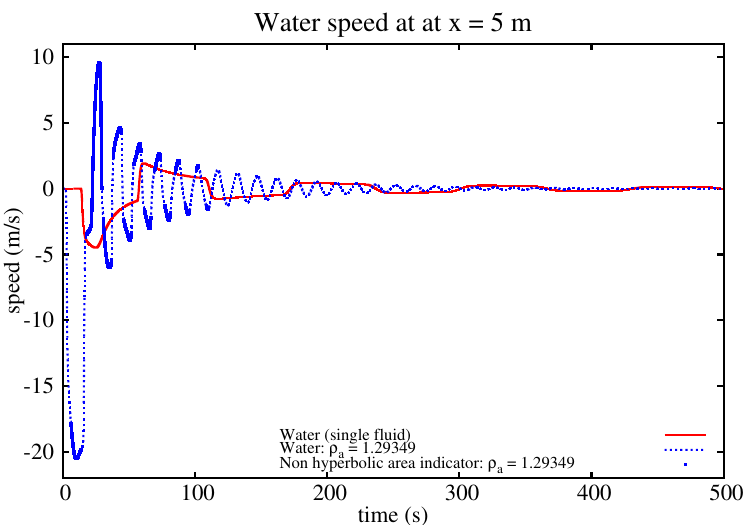}
} 

\subfigure[Air speed \label{Bar5AirSpeed}]
{
\includegraphics[height = 5.8cm]{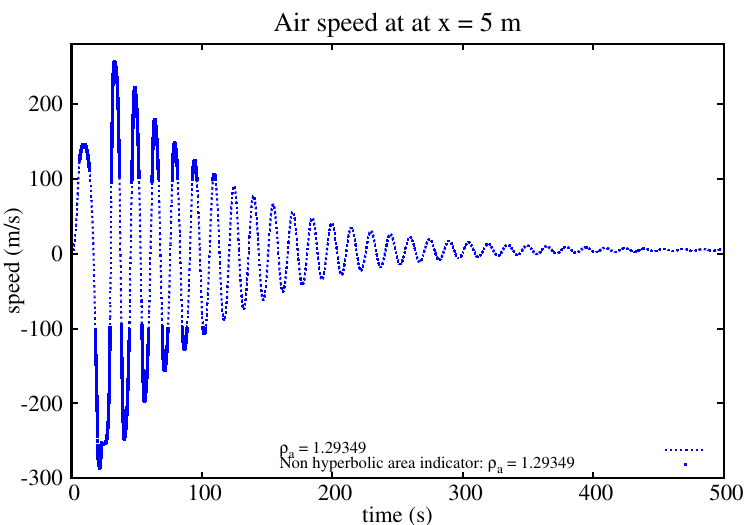}
}
\caption{Non constant downstream water discharge: numerical results at $x=5\: m$\label{Bar5}}
\end{figure}

A second numerical tests is done in order to confirm the previous conclusion only for  
 $\rho_{a} = 1.29349\:  kg/m^3$ fixed at upstream and downstream end (since the other cases do not show real differences) that we compare to the 
single fluid in the following settings:
the horizontal pipe of circular cross-section of $2 \: m^2$ is $100\: m $ long. 
The altitude of the upstream end of the pipe is $0 \: m$.  
We start the numerical simulation with a  steady state with the water height $h_{init}(x) = 0.05\:  m$ 
(i.e. $A_{init}(x) = 2.0\, 10^{-2}\: m^2$), $M_{init}(x):=\frac{\rho_{a}}{\rho_0}(S(x)-A_{init}(x)) = 4.00\,  10^{-3}\: m^2$, 
the water discharge $Q_{init}(x) = 0\: m^3/s$, the air discharge $D_{init}(x) = 0\: m^3/s$, for all $x\in [0,100]$.   
The upstream water discharge is keep constant to $Q_0(t) = 0\: m^3/s$ and the downstream water discharge is set to 
$Q_{N+1}(t)=  - 5 t \mathds{1}_{[0,20]}(t)$. 
Other parameters are: 
$$
\begin{array}{lcl}
\hbox{Discretisation points}& : & 50,\\ 
\hbox{Delta x }(m)& : &0.5,\\
\hbox{CFL }& : & 0.95,\\
\hbox{Simulation time } (s) & : &500.
\end{array}
$$
On \resim\ref{Bar5}, we represent the water height, the water and air discharge.  
The same phenomea are reproduced also in the present numerical experiment, i.e. acceleration phenomenon and oscillations are also 
present.  Thus, this confirm the previous observations about such phenomena.

\subsection{Numerical order of the discretization error of the kinetic numerical scheme}\label{ordre}
In the previous numerical experiments, we have shown through several examples that the scheme even if the system is not hyperbolic seems to be 
stable. We have also pointed out that the presence of oscillations was not due to the loss of hyperbolicity. 
Now, we will check the sensitivity of the two-layer kinetic scheme with respect to decreasing spatial mesh size  and show that no 
geometrically (as we may also emphasize in the previous numerical tests) growing instabilities are observed.

In order to obtain a qualitative behavior of the scheme and to compute a ``numerical'' 
order of the discretization error of the two-layer kinetic numerical scheme, 
we present the ``dam-break'' problem which is defined as initial state  
by for all $x \in[0,L]$, $$\,\,\,Q_{init}(x) = 0\: m^3/s,\, {\rm and} \quad h_{init}(x) = 
1 \mathds{1}_{[0,L/2)}(x)+ 0.5 \mathds{1}_{(L/2,L]}(x) \: m$$
for the water layer and $D_{init}(x) =0\: m^3/s$ and $M_{init}(x)=\rho_{init}(x)/\rho_0 (S-A_{init}(x))$ for the air layer. 
We prescribe at the upstream and the downstream end a null water discharge while boundary conditions will be provided by Equation 
\eqref{BCair} to compute $D_0^n$ and $D_{N+1}^n$ for each time $t_n$.

\begin{figure}[H]
\subfigure[Initial state  \label{NOTempsT0IS}]{\includegraphics[height = 7.8cm]{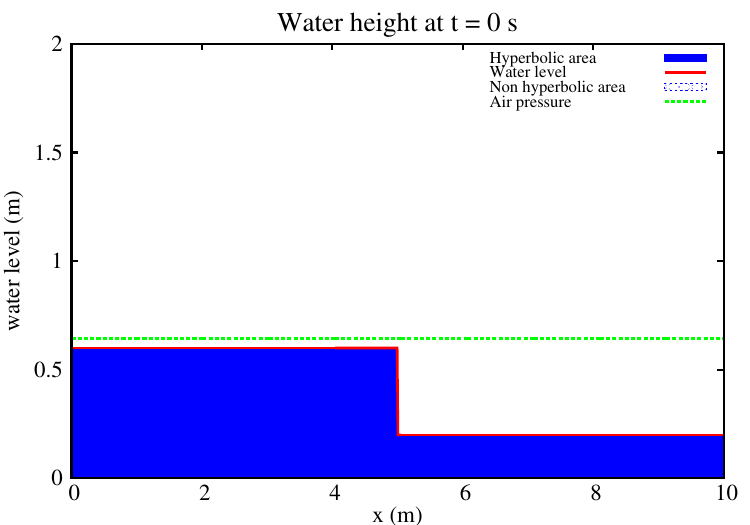}}
\subfigure[$L_t^1(L^2_x)$ norm of the water level \label{NOTempsT0L2xtHeight}]{\includegraphics[height = 7.8cm]{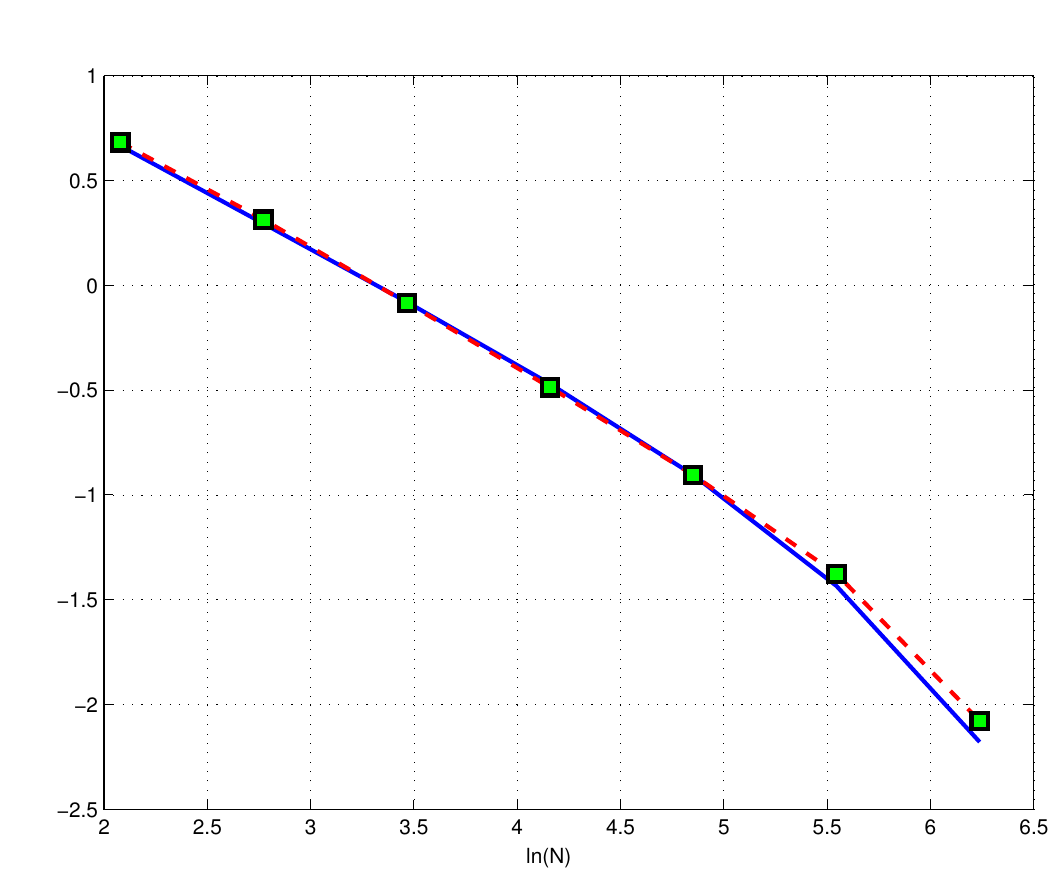}}
\caption{Initial state for the ``dam-break'' and the $L_t^1(L^2_x)$ norm of the error (in $\log$ scale). }\label{NOTempsT0}
\end{figure}

The horizontal pipe is a circular one of diameter $2 \:m$ and  $L = 10 \: m$ long. 
The altitude of the upstream end of the pipe is $0 \: m$.
The upstream and downstream discharge is kept constant equal to $0\: m^3/s$.  
The air density is chosen such as the two-layer system is mainly non  hyperbolic as we can see on 
\resim\ref{NOTempsT0IS}: 
it will be set to $\rho_{a} = 1.29349\:  kg/m^3$.
We have computed the ``exact'' numerical flow for all time from $t=0 \:s$ until $t= 2 \: s$ 
with a uniform  discretization of $ N = 1024$ mesh points. 
We have then computed the $L^2$ norm of the difference between the water level  computed by the numerical 
kinetic scheme for different value of the mesh points and the ``exact'' numerical solution 
at time $t = 0.02  \:s$, $t = 0.2  \:s$ and $t= 1.688 \:s$.  Results are then compared to the ``single fluid'' case.

Other parameters are: 
$$
\begin{array}{lcl}
\hbox{CFL }& : & 0.95\,,\\
\hbox{Simulation time } (s) & : &2 \,,\\
\hbox{Delta t }(s)& : & 10^{-5} .
\end{array}
$$
We have fixed the time step in order to obtain the solution of the ``single fluid'' and two layer model at the same time. Thus, we ensure that 
$\Delta t = 10^{-5} $ satisfies the condition $10^{-5} >  \frac{CFL \; \Delta x}{\max_{\alpha,i=0,\ldots,N+1}\big(| u_{\alpha,i}^n| +\sqrt{3}b_{\alpha,i}^n\big)}$.

\begin{figure}[H]
\subfigure[Two-layer model  \label{NOTempsT1AIR}]{\includegraphics[height = 5.8cm]{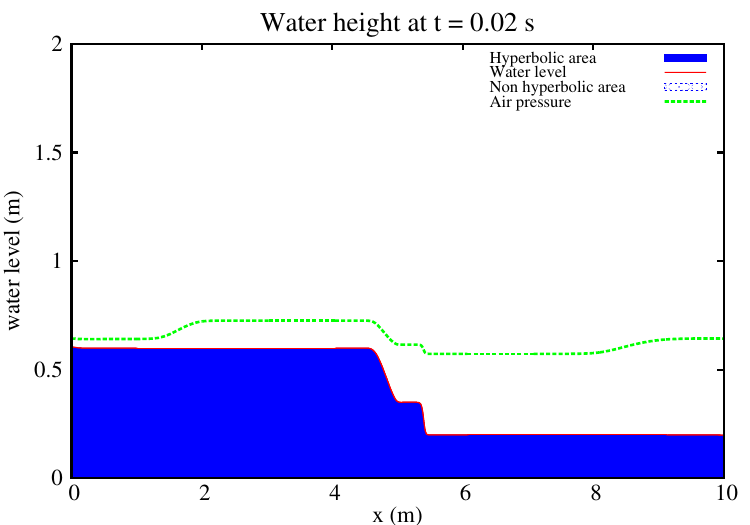}}\\
\subfigure[The single fluid \label{NOTempsT1WATER}]{\includegraphics[height = 5.8cm]{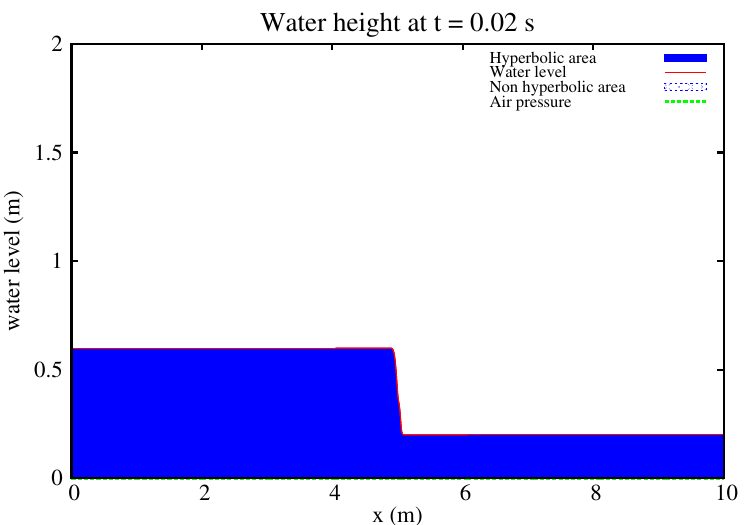}}\\
\subfigure[Numerical order  \label{NOTempsT1ORDER}]{\includegraphics[height = 5.8cm]{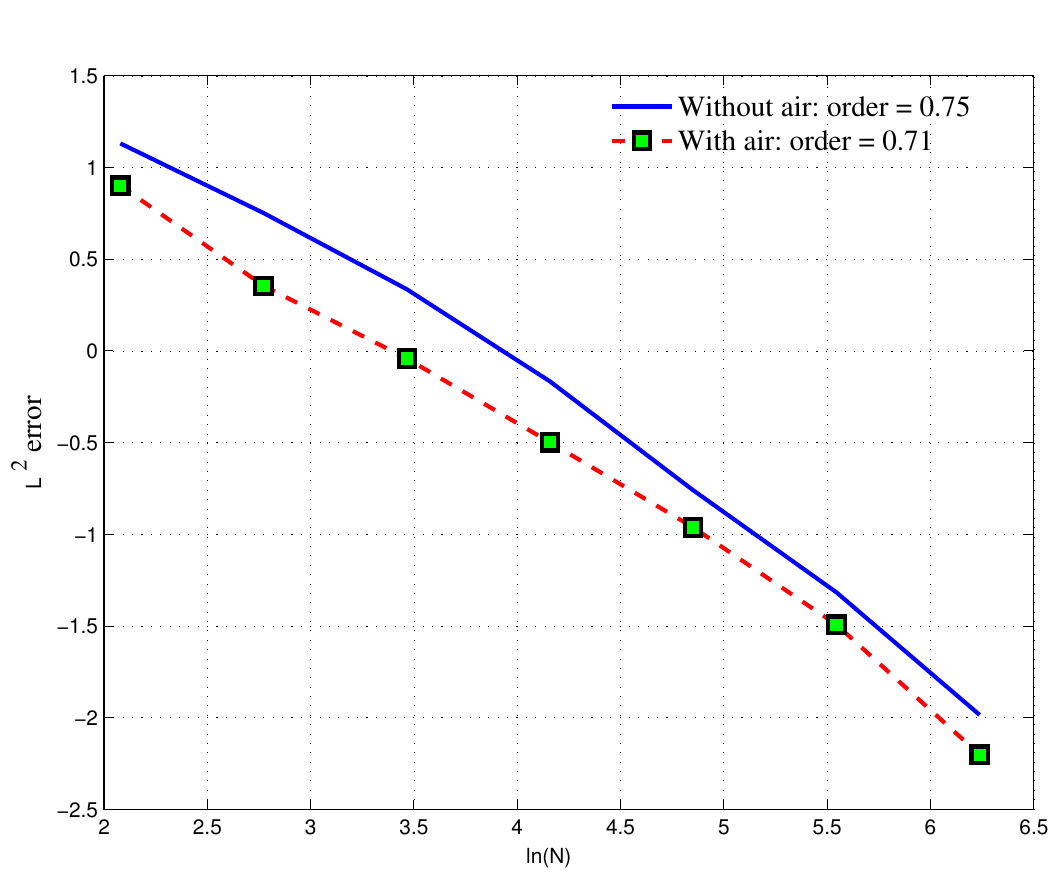}}
\caption{The ``dam-break'' problem at $t=0.02\; s$. $L^2$ norm of the error  (in log  scale).}\label{NOTempsT1}
\end{figure}

\begin{figure}[H]
\subfigure[Two-layer model  \label{NOTempsT2AIR}]{\includegraphics[height = 5.8cm]{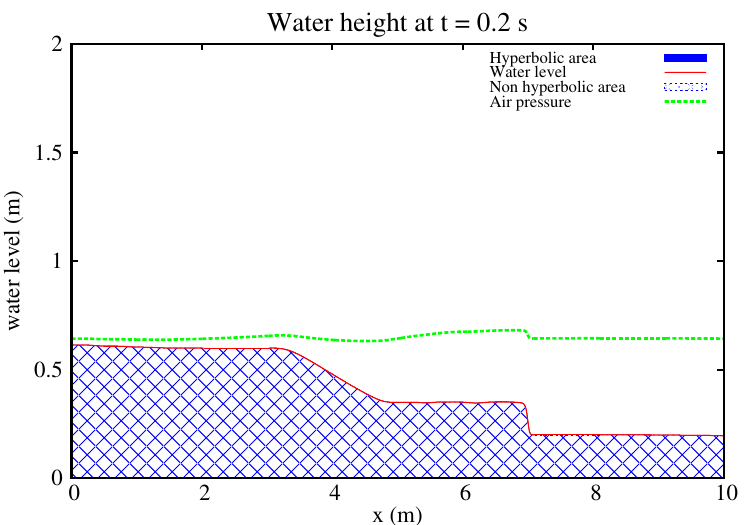}}\\
\subfigure[The single fluid \label{NOTempsT2WATER}]{\includegraphics[height = 5.8cm]{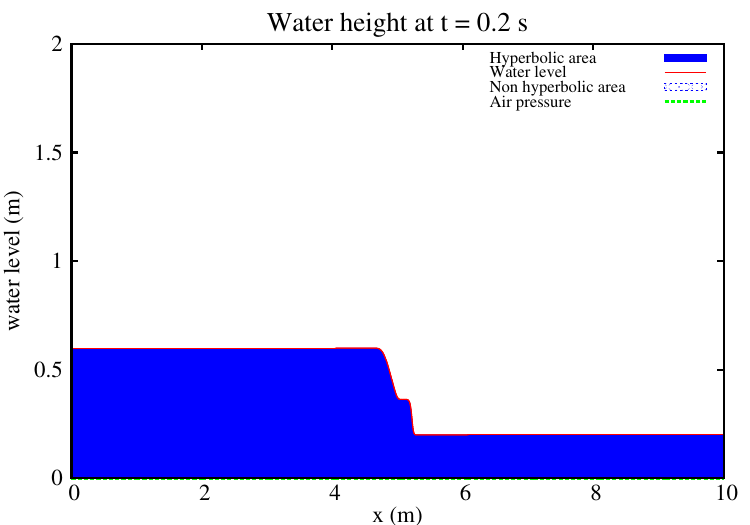}}\\
\subfigure[Numerical order  \label{NOTempsT2ORDER}]{\includegraphics[height = 5.8cm]{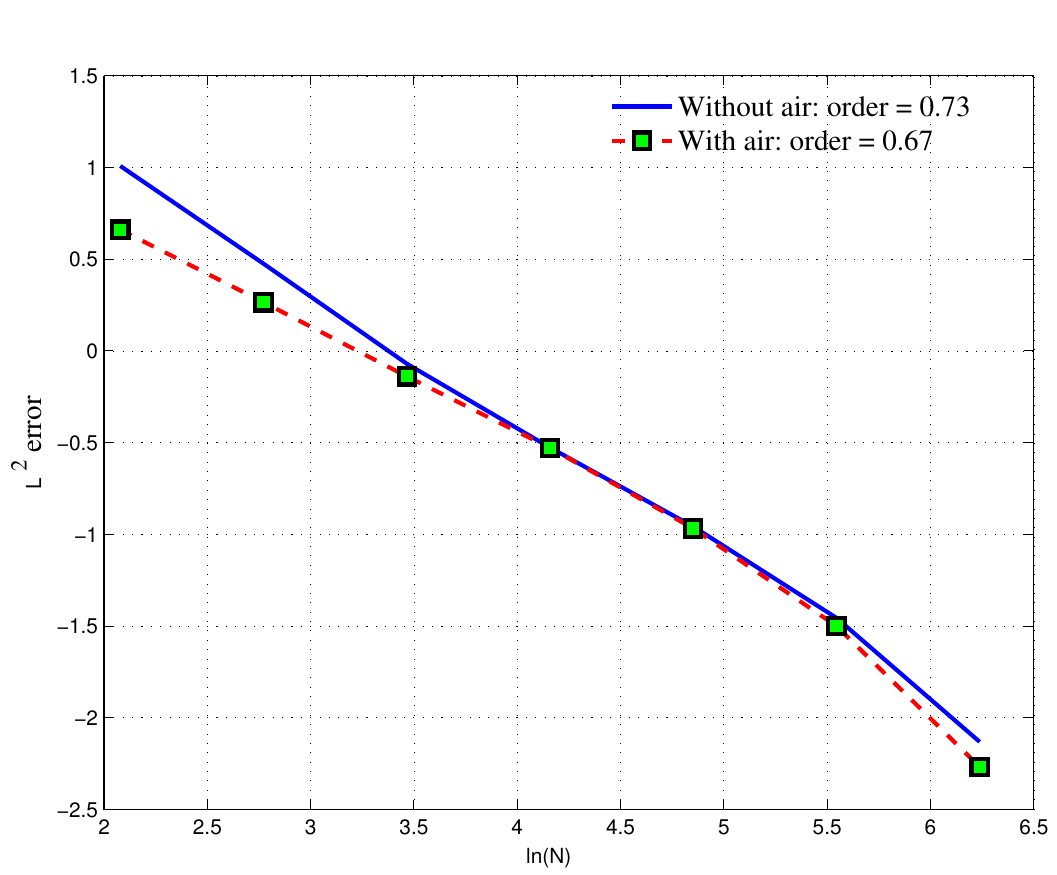}}
\caption{The ``dam-break'' problem at $t=0.1\; s$. $L^2$ norm of the error  (in log scale).}\label{NOTempsT2}
\end{figure}

\begin{figure}[H]
\subfigure[Two-layer model  \label{NOTempsT3AIR}]{\includegraphics[height = 5.8cm]{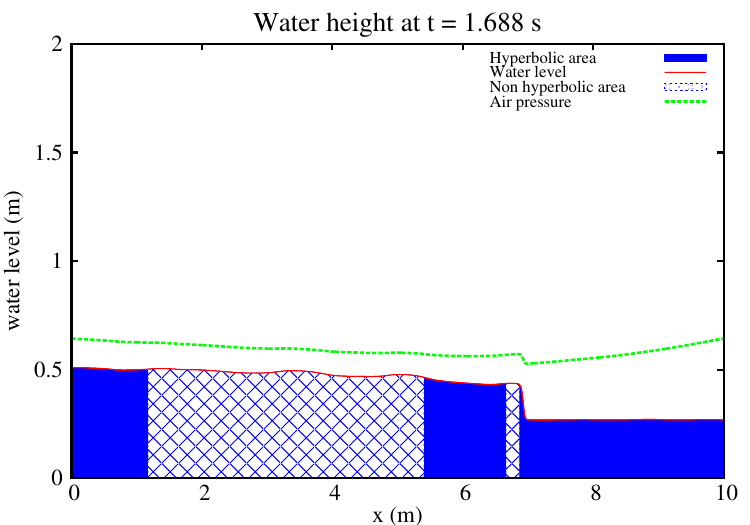}}\\
\subfigure[The single fluid \label{NOTempsT3WATER}]{\includegraphics[height = 5.8cm]{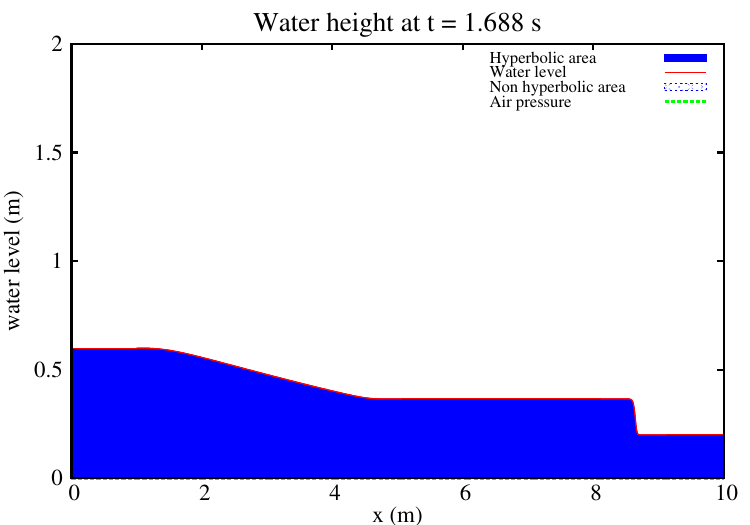}}\\
\subfigure[Numerical order  \label{NOTempsT3ORDER}]{\includegraphics[height = 5.8cm]{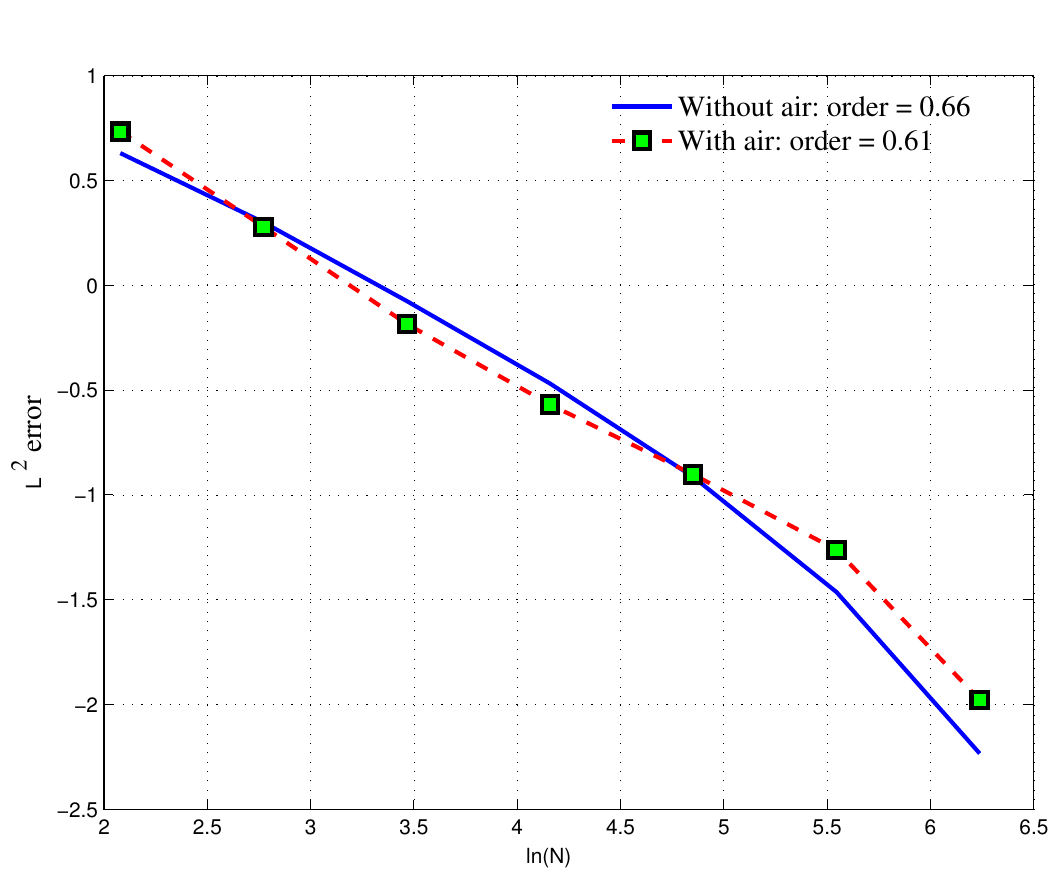}}
\caption{The ``dam-break'' problem at $t=1.688\; s$. $L^2$ norm of the error  (in log  scale).}\label{NOTempsT3}
\end{figure}

For each time $t=0.02$ (see \resim\ref{NOTempsT1}), $t=0.2$ (see \resim\ref{NOTempsT2}) and $t=1.688$ (see \resim\ref{NOTempsT3}), 
we display the shape of the water level for both cases and the obtained $L^2$ numerical order. \resim\ref{NOTempsT1} displays the case where 
the system is hyperbolic, \resim\ref{NOTempsT2} displays the case where the system is non hyperbolic and finally on \resim\ref{NOTempsT3}, the 
system is partially hyperbolic.

Even if the system loses its hyperbolicity, the numerical order of the two-layer kinetic scheme is very 
close to the case of ``single fluid'' for 
the $L^2$ norm as well as for the $L_t^1(L^2_x)$ (see \resim\ref{NOTempsT0L2xtHeight}) which provides a 
numerical stability as the space step goes to $0$. 

Finally, even if we do not know if the presented numerical scheme 
preserves the entropy inequalities, this  stability should come from the fact that 
the system at continuous level preserves an entropy equality (see Theorem \ref{BilayerModelThm}).

\section{Conclusion}
This study is of course a first step in the comprehension and the modelisation of the air role in the transient flows in closed pipes. 
The kinetic scheme seems to be well-adpated to treat a two layer model even if the obtained partial differential system is conditionally hyperbolic.

We have now to treat the air entrapment pocket 
which may be encountered in a rapid filling process of a closed pipe: in this case a portion of the pipe will be completely filled and an air pocket will be entrapped. 
This is the next step of our research work, since in previous works  we have derived a model for pressurised flow or mixed flows in closed pipes
 {\it without  taking into account the role of the air} and proposed a Roe like Finite Volume method and a kinetic scheme 
 \cite{BEG09_1,BEG09_2, BEG11_1, BEG11_2,BG07}.

Next, we have to deal with the evaporation/condensation of gas or water when it is not supposed to be isothermal and 
later, because in the water hammer phenomenon, large depression may occur, we have to deal with the natural cavitation problem.

\section*{Acknwoledgements}
This work is supported by the ``Agence Nationale de la Recherche'' referenced by  ANR-08-BLAN-0301-01 and 
the  second author was  supported by  the ERC Advanced Grant FP7-246775 NUMERIWAVES. This work was finalized while the third author was visiting 
BCAM--Basque Center for Applied Mathematics, Derio, Spain, and partially supported by  the ERC Advanced Grant FP7-246775 NUMERIWAVES. The third author wishes to thank Enrique Zuazua for his kind hospitality.

The authors thank the referees for their valuable remarks which led to substantial improvement of the first version of this paper.

\bibliographystyle{plain}


\end{document}